\documentclass[10pt]{amsart}

\usepackage{graphicx}
\usepackage{mathtools}
\usepackage{amsmath} 
\usepackage{amsfonts}
\usepackage{amssymb}
\usepackage{amsthm}
\usepackage{latexsym}
\usepackage{mathtools}
\usepackage{enumerate}
\usepackage{mathrsfs}
\usepackage{color}

\newtheorem{theorem}{Theorem}[section]
\newtheorem{lemma}[theorem]{Lemma}
\newtheorem{definition}[theorem]{Definition}

\newtheorem{remark}[theorem]{Remark}

\numberwithin{equation}{section}

\newtheorem*{theorem*}{Theorem}
\newtheorem{prop}[theorem]{Proposition}

\newtheorem*{bartnikMIN*}{Bartnik mass of minimal Bartnik data}
\newtheorem*{bartnikMINAH*}{AH Bartnik mass estimate of minimal Bartnik data}
\newtheorem*{bartnikMINC*}{Charged Bartnik mass of minimal Bartnik data}

\newcommand{\pr}{\partial}
\newcommand{\veps}{\varepsilon}

\newcommand{\ovr}[1]{\overline{#1}}

\usepackage{enumerate}
\usepackage{mathrsfs}

\usepackage{color}
\usepackage{tikz}
\usetikzlibrary{matrix}
\usepackage{pgfplots}

\usetikzlibrary{intersections}
\usetikzlibrary{shapes,arrows}
\usetikzlibrary{calc}
\usetikzlibrary{shapes.geometric}
\usetikzlibrary{shapes.arrows}
\usetikzlibrary{arrows.meta}
\usetikzlibrary{decorations.markings}
\usetikzlibrary{fit}
\usetikzlibrary{hobby}
\usetikzlibrary{patterns}
\usetikzlibrary{shapes.misc}
\pgfplotsset{compat=1.12}
\usepgfplotslibrary{patchplots}

\def\II{\textnormal{I \hspace{-.4em}I}}
\def\tr{\textnormal{tr}}
\def\dv{\textnormal{div}}
\def\M{\mathscr{M}}
\def\Lap{\Delta}
\def\Ric{\textnormal{Ric}}
\def\diver{\textnormal{div}}
\def\R{\mathbb{R}}

\def\R{\mathbb{R}}

\def\S{\Sigma}
\def\({\left(}
\def\){\right)}
\def\a{\alpha}
\def\b{\beta}
\def\to{\longrightarrow}
\def\Hb{\mathbb{H}}
\def\s{\sigma}
\def\w{\omega}
\def\g{\gamma}
\def\ADM{\textnormal{ADM}}
\def\AH{\textnormal{AH}}
\def\CH{\textnormal{CH}}
\def\bS{\mathbb{S}}
\def\m{\mathfrak{m}}
\def\AdS{\textnormal{AdS}}
\newcommand{\definedas}{\mathrel{\raise.095ex\hbox{\rm :}\mkern-5.2mu=}}
 \newcommand{\asdefined}{\mathrel{=\mkern-5.2mu\raise.095ex\hbox{\rm :}}}
\def\defeq{\definedas}

\begin{document}
\title[Extensions and Bartnik mass estimates]{A survey on extensions of Riemannian manifolds and Bartnik mass estimates}

\author[{Cabrera Pacheco}]{Armando J. {Cabrera Pacheco}}
\address{Department of Mathematics, Universit\"at T\"ubingen, T\"ubingen 72076, Germany}
\curraddr{}
\email{cabrera@math.uni-tuebingen.de}
\thanks{}

\author[Cederbaum]{Carla Cederbaum}
\address{Department of Mathematics, Universit\"at T\"ubingen, T\"ubingen 72076, Germany}
\curraddr{}
\email{cederbaum@math.uni-tuebingen.de}
\thanks{}


\date{}

\begin{abstract}
Mantoulidis and Schoen developed a novel technique to handcraft asymptotically flat extensions of Riemannian manifolds $(\S \cong \bS^2,g)$, with $g$ satisfying $\lambda_1 \defeq \lambda_1(-\Lap_g + K(g))>0$, where $\lambda_1$ is the first eigenvalue of the operator $-\Lap_g+K(g)$ and $K(g)$ is the Gaussian curvature of $g$, with control on the ADM mass of the extension. Remarkably, this procedure allowed them to compute the Bartnik mass in this so-called minimal case; the Bartnik mass is a notion of quasi-local mass in General Relativity which is very challenging to compute. In this survey, we describe the Mantoulidis--Schoen construction, its impact and influence in subsequent research related to Bartnik mass estimates when the minimality assumption is dropped, and its adaptation to other settings of interest in General Relativity.
\end{abstract}

\phantom{a}
\vspace{2cm}
\maketitle

\section{Mathematical Relativity} \label{sec-intro}

In the context of Mathematical Relativity, a Riemannian manifold $(M^{3},\g)$ with non-negative scalar curvature $R(\g)$ represents a ``time-symmetric time-slice'' of a spacetime satisfying the ``dominant energy condition''. From the point of view of the Cauchy problem in General Relativity, we can consider $(M^{3},\g)$ as a time-symmetric \emph{initial data set} for the Einstein Equations, and the dominant energy condition $R(\g) \geq 0$ as a compatibility condition that will allow this manifold to evolve into a physically relevant spacetime containing $(M^{3},\g)$ as a submanifold. The initial value problem in General Relativity has been widely studied and we refer the interested reader to the books by Choquet-Bruhat~\cite{CB-book} and Ringstr\"om~\cite{Ringstrom}.

In this survey, we are mostly interested in \emph{asymptotically flat Riemannian manifolds}. Consider a gravitating system in which all the matter is contained inside a compact region $\ovr{K}$. Then, as we move away from $\ovr{K}$, the gravitational field will get weaker, and if we go far enough away, we will increasingly feel as if there was no gravitation exerted on us anymore. For a time-symmetric initial data set $(M^{3},\g)$ modeling this gravitating system, this means that as we move away from $\ovr{K}\subset M^{3}$, the metric $\g$ approaches the Euclidean metric $\delta$. More precisely, a Riemannian manifold $(M^3,\g)$ is said to be \emph{asymptotically flat} if $M^{3} \setminus \ovr{K}$ is diffeomorphic to $\R^3 \setminus \ovr{B_1(0)}$ via a map $\Phi$ which defines a set of \emph{coordinates at infinity}, and in these coordinates the metric takes the form $\g_{ij}=\delta_{ij}+\mathcal{O}_{2}(|x|^{-p})$ with $p>\frac{1}{2}$. In addition, we require the scalar curvature $R(\g)$ of $\g$ to be integrable over $M^3$.

For example, consider the famous \emph{Schwarzschild spacetime of mass $m>0$}, which models the gravitational field surrounding a static, non-rotating, rotationally symmetric black hole or celestial body of mass $m$ in vacuum. The spacetime manifold can be written as $\ovr{M}^4=\R \times M^3=\R \times (2m,\infty) \times \bS^2 \ni(t,r,x)$, with the Lorentzian (spacetime) metric defined by 
\begin{align*}
\ovr{\g}_m = -\(1 - \frac{2m}{r}\) dt^2 +  \(1 - \frac{2m}{r}\)^{-1} dr^2 + r^2 g_*,
\end{align*}
where $g_*$ denotes the standard round metric on $\bS^2$. The Riemannian manifold $M^3$ defined by $\{ t=0 \}$ is called the \emph{(spatial) Schwarzschild manifold of mass $m$}, and clearly, it consists of $M^3=(2m,\infty) \times \bS^2$ with the Riemannian metric
\begin{align}\label{eq-Schwarzschild-metric}
\g_m =  \(1 - \frac{2m}{r}\)^{-1} dr^2 + r^2 g_*.
\end{align}
It can be checked that $R(\g_m) \equiv 0$, and evidently $\g_m$ approaches the Euclidean metric $\delta$ as $r \to \infty$.

Remarkably, by imposing the decay conditions on the metric described above, an asymptotically flat manifold $(M^3,\g)$ behaves in a similar way to a Schwarzschild manifold and we can ``detect'' the total mass of the gravitating system it models --- or its \emph{total mass} for short --- from its asymptotic behavior. This notion of total mass is called ADM mass (defined by Arnowitt, Deser, and Misner in~\cite{ADM}) and it is denoted by $m_{\ADM}(M^3,\g)$; whenever there is no risk of confusion we will drop the reference to the underlying manifold $M^{3}$. As the reader will discover, for this survey we do not need to use the explicit definition of the ADM mass, rather we will only need to keep in mind that it is measured from the asymptotic behavior of $(M^{3},\g)$ and that a Schwarzschild manifold of mass $m$ has, in fact, $m_{\ADM}(\g_m)=m$. We avoid any digress towards technical definitions and convergence issues by referring the interested reader to~\cite{ADM}, and to the works of Bartnik~\cite{Bartnik-ADM} and Chru\'sciel~\cite{Chrusciel}.

With a notion of total mass at hand, the next natural step is to determine if, under some reasonable physical assumptions, it is non-negative. This is known as the Positive Mass Theorem and was proven independently by Schoen and Yau~\cite{SY1} using minimal surface theory and by Witten~\cite{Witten} using spinors.

\begin{theorem*}(Positive Mass Theorem) \label{thm-PMT}
Let $(M^3,\g)$ be a complete, asymptotically flat manifold with $R(\g) \geq 0$. Then $m_{\ADM}(M^3,\g) \geq 0$, and equality holds if and only if $(M^3,\g)$ is isometric to the Euclidean space $(\R^3,\delta)$.
\end{theorem*}

The important Riemannian Penrose Inequality states that the mass of an asymptotically flat Riemannian manifold is bounded below in terms of the area of $2$-surfaces representing black holes. In the time-symmetric case (that is, for Riemannian manifolds), black holes are represented by compact minimal surfaces, i.e., compact surfaces with vanishing mean curvature, which are called \emph{horizons}. A horizon is said to be \emph{outer-minimizing} if it minimizes area among all surfaces that enclose it. As an example, consider the surface $\{ r=2m \}$ in a Schwarzschild manifold of mass $m>0$ --- but note that, as it is written, $\{r=2m\}$ is in fact not part of the Schwarzschild manifold, and if we naively try to include it as an inner boundary, $\g_m$ would appear to degenerate there. Nevertheless, one can perform a change of coordinates that allows a smooth extension of $\g_m$ to the boundary (see Section~\ref{subsec-gluing}).

\begin{theorem*}[Riemannian Penrose Inequality] \label{thm-RPI}
Let $(M^3,\g)$ be an asymptotically flat manifold with $R(\g) \geq 0$ and such that its boundary $\pr M$ is an outer-minimizing horizon. Then
\begin{align}\label{eq-RPI-bound}
m_{\ADM}(M^3,\g) \geq \sqrt{\frac{|\pr M|_{\g}}{16 \pi}},
\end{align}
where $|\pr M|_{\g}$ represents the area of the boundary of $M$ with respect to $\g$. Moreover, equality holds if and only if $(M^3,\g)$ is isometric to a Schwarzschild manifold of positive mass.
\end{theorem*}

The Riemannian Penrose Inequality has been proven by Huisken and Ilmanen~\cite{HI} using a weak version of inverse mean curvature flow (IMCF) when the horizon has a single connected component, and by Bray~\cite{Bray} using conformal flow and thereby allowing for multiple connected components. There is plenty of literature available about the Riemannian Penrose Inequality and the more general (spacetime) Penrose Inequality; we particularly recommend to consult Mars' survey article~\cite{Mars} which also covers other geometric inequalities from General Relativity that will appear in this survey.

An interesting question in Mathematical Relativity is whether it is possible to define a local notion of mass contained in a small region $\Omega \subset M^3$. For various conceptual as well as physical reasons, it turns out to be more convenient in many cases to work with localized notions of mass defined in terms of the geometry of $\S \defeq \pr \Omega$. Such notions of mass are referred to as \emph{quasi-local mass} notions. We now define a highly important example of such a notion.

\begin{definition}
Let $\S$ be a surface in a Riemannian manifold $(M^{3},\g)$ with induced metric $g$ and mean curvature $H$. The \emph{Hawking mass} of the triplet $(\S,g,H)$ is defined as
\begin{align*}
\m_\textnormal{H}(\S,g,H) \defeq \sqrt{\frac{|\S|_g}{16 \pi}}\( 1 - \frac{1}{16 \pi} \int_{\S} H^2 \, d\s \),
\end{align*}
where $|\S|_g$ is the area of $\S$ with respect to the metric $g$ and $d\s$ denotes the area form on $\S$.
\end{definition}

Note that, using the notion of Hawking mass, we can simply rewrite \eqref{eq-RPI-bound} as
\begin{align} \label{eq-Hawking-bound}
m_{\ADM}(M^3,\g) \geq \m_\textnormal{H}(\pr M, \g\vert_{\pr M}, H\equiv 0).
\end{align}

Another important notion of quasi-local mass, which will be one of the main subjects in this survey, is the \emph{Bartnik mass}, defined by Bartnik in~\cite{Bartnik-mass}. 

\begin{definition}
A triple $(\S \cong \bS^2,g,H)$, where $g$ is a Riemannian metric and $H \geq 0$ is a smooth function on $\S$ is called \emph{Bartnik data}. Bartnik data $(\S \cong \bS^2,g,H)$ with $H \equiv 0$ are called \emph{minimal Bartnik data}.
\end{definition}

For given Bartnik data $(\S \cong \bS^2,g,H)$, we consider the set of admissible extensions $\mathcal{A}$, consisting of all asymptotically flat manifolds $(M^3,\g)$ such that $R(\g) \geq 0$, the boundary $\pr M$ is outer-minimizing, and $(\pr M,\g\vert_{\pr M} )\cong (\S,g)$ with mean curvature $H$. 

\begin{definition} \label{def-Bartnik}
The \emph{Bartnik mass} of Bartnik data $(\S \cong \bS^2,g,H)$ is defined as
\begin{equation} \label{eq-Bartnik-mass}
\m_\textnormal{B} (\S \cong \bS^2,g,H) \defeq \inf\{ m_{\ADM}(M^3,\g)\,\vert\,{(M^3,\g) \in \mathcal{A}}  \}.
\end{equation}
\end{definition}

Evidently, by definition, the Bartnik mass is extremely complicated to compute. However, notice that when $H=0$, an admissible extension satisfies the conditions of the Riemannian Penrose Inequality (Theorem~\ref{thm-RPI}), and hence, the Hawking mass provides a \emph{lower bound} for the Bartnik mass. The original derivation leading to the notion of Bartnik mass does in fact not readily extend to the case of minimal Bartnik data. We suppress this delicate debate here.

By Huisken and Ilmanen's proof of the Riemannian Penrose Inequality~\cite{HI}, we know that
\begin{align}\label{genPI}
\m_\textnormal{H}(\S \cong \bS^2,g,H) \leq m_{{\ADM}}(M^{3},\g)
\end{align}
for any admissible extension $(M^{3},\g)$ of given Bartnik data $(\S \cong \bS^2,g,H)$, again with equality if and only if the Bartnik data $(\S \cong \bS^2,g,H)$ suitably embed into a Schwarzschild manifold of positive mass. This allows to compute the Bartnik mass in the case of spherical symmetry --- i.e., in the case of Bartnik data that isometrically embed as a centered round sphere into a Schwarzschild manifold of positive mass $m$; in that case, the Bartnik and Hawking mass both coincide with the mass parameter $m$.

From Definition \ref{def-Bartnik}, it is apparent that computing the Bartnik mass of given Bartnik data is a challenging task. However, also from its definition, we see that any extension of given Bartnik data will provide an \emph{upper bound} for the Bartnik mass. This in part motivated Bartnik himself to study asymptotically flat ``quasi-spherical'' extensions of intrinsically round Riemannian $2$-spheres $(\S,g)$ in~\cite{Bartnik-QS}. Other extensions using similar techniques to those in~\cite{Bartnik-QS} are the works of Smith and Weinstein \cite{SW1,SW2}, Shi and Tam \cite{ShiTam02}, Smith \cite{Smith} and Lin \cite{Lin}. Building on the work of Lin in \cite{Lin}, Lin and Sormani constructed asymptotically flat extensions to estimate the Bartnik mass in \cite{LS}.

In 2016, Mantoulidis and Schoen~\cite{MS} first succeeded in computing the Bartnik mass of non-spherically symmetric Bartnik data. More concretely, they were able to compute the Bartnik mass of minimal Bartnik data $(\S \cong \bS^2,g,H \equiv 0)$ in the case when $g$ satisfies $\lambda_1 \defeq \lambda_1(-\Lap_g + K(g))>0$, where $\lambda_1$ denotes the first eigenvalue of the operator $-\Lap_g + K(g)$, and $K(g)$ is the Gaussian curvature of $g$. To do so, they constructed asymptotically flat extensions of minimal Bartnik data $(\S \cong \bS^2,g,H \equiv 0)$ in such a way that they can control their ADM masses, ultimately leading them to prove that  $\m_\textnormal{B}(\S \cong \bS^2,g,H \equiv 0) = \m_\textnormal{H}(\S \cong \bS^2,g,H \equiv 0) $. Their work has inspired related constructions, among other results, that provide estimates for Bartnik data with $H \neq 0$. The Mantoulidis--Schoen construction together with its adaptations and modifications, specifically those useful to obtain Bartnik mass estimates, will be the main topic of this survey.

Very recently, building on work by Shi and Tam  \cite{ShiTam02} and by Miao~\cite{Miao-locRPI}, Miao and Xie constructed extensions of vanishing scalar curvature for Bartnik data of positive Gaussian curvature in~\cite{MXie}. In particular, for the minimal case they recover that $\m_\textnormal{B}(\S \cong \bS^2,g,H \equiv 0) = \m_\textnormal{H}(\S \cong \bS^2,g,H \equiv 0)$. See \cite{MXie} for details.

For the sake of completeness, let us note that it is common in the literature, for very good reasons, to request that admissible extensions for given Bartnik data satisfy that $(M^3,\g)$ contains no other minimal surfaces (homologous to $\pr M$), except perhaps for $\pr M$ --- instead of requesting the outer-minimizing condition on the boundary that we used above, see e.g.~\cite{Bray,HI}, and the work of Bray and Chru\'sciel~\cite{BC}. In this survey, it will not make any difference which condition is chosen since all the extensions considered here will satisfy both. 

We further remark that there are in fact also other definitions of Bartnik mass in the literature. This raises natural questions about the relations or equivalences among them, which will not be discussed here. We refer the interested reader to the works of Jauregui~\cite{Jauregui-B} and McCormick~\cite{McCormick-B}, which specifically address these questions of equivalence of various definitions of Bartnik mass. For more information on Bartnik mass, see also the review~\cite{GMS} by Galloway, Miao, and Schoen.

\subsection*{Acknowledgments} We would like to thank the anonymous referee, Jeffrey  L.~Jauregui, Christos Mantoulidis, Stephen McCormick, Pengzi Miao, and Richard M.~Schoen for helpful comments and interesting conversations regarding this survey. We thank Axel Fehrenbach for creating the pictures for this survey.

The authors are grateful to the Carl Zeiss Foundation for its generous support. Work of CC is supported by the Institutional Strategy of the University of T\"ubingen (Deutsche Forschungsgemeinschaft, ZUK 63). \medskip

The structure of this survey is as follows. In Section~\ref{sec-MS}, we intuitively explain the Mantoulidis--Schoen construction and describe works related to Bartnik mass estimates in a historical order. In Section~\ref{section-gluingmethods}, we describe in detail key steps used to construct these types of extensions and how to obtain Bartnik mass estimates from them. In Section~\ref{sec-Bartnik-est}, we describe in more detail the specific Bartnik mass estimates obtained with this method. Finally, in Section~\ref{sec-conclusions}, we discuss the current state of the art and other related problems.

\subsection*{Some useful considerations}
In what follows, we will give a brief description of higher dimensional Riemannian manifolds in the context of Mathematical Relativity, of asymptotically hyperbolic Riemannian manifolds, and of charged Riemannian manifolds arising as time-symmetric initial data sets for the Einstein--Maxwell Equations. We recommend that the reader skips this part upon first reading and returns to it as it becomes necessary.

\subsection{Higher dimensions} \label{subsec-intro-higherdim}

A Riemannian manifold $(M^{n+1},\g)$ of dimension $n+1$ is said to be \emph{asymptotically flat} if, as before, outside a compact set, it is diffeomorphic to $\R^{n+1}\setminus\ovr{B_1(0)}$. In this coordinate chart at infinity, we ask that $\g_{ij}=\delta_{ij}+\mathcal{O}_{2}(|x|^{-p})$ for $p > \frac{n-2}{2}$. In addition, $R(\g)$ must be integrable. As in the $3$-dimensional case, by imposing these decay conditions on the metric, the total mass $m_{\ADM}(M^{n+1},\g)$ is well-defined. The Positive Mass Theorem also holds in higher dimensions. It was proven by Schoen and Yau for $2 \leq n < 7$~\cite{SY2} (the case $n=7$ follows from the work of Smale~\cite{SN}) and for spin manifolds in general dimension by Witten~\cite{Witten}. Recently, Schoen and Yau established the Positive Mass Theorem in all dimensions~\cite{SY3}. The Riemannian Penrose Inequality in higher dimensions, proven by Bray and Lee~\cite{BL}, states that for an asymptotically flat manifold $(M^{n+1},\g)$ with $R(\g) \geq 0$ (for $2 \leq n < 7$), such that its boundary $\pr M$ is an outer-minimizing horizon, we have
\begin{align} \label{eq-RPI-hd}
m_{\ADM}(M^{n+1},\g) \geq \frac{1}{2}  \(  \frac{|\pr M|_{\g}}{ w_{n}} \)^{\frac{n-1}{n}},
\end{align}
where $\w_n$ denotes the volume of the standard sphere $\bS^n$ and $|\pr M|_{\g}$ represents the area of the boundary of $M$ with respect to $\g$ as before. Equality holds if and only if $(M^{n+1},\g)$ is isometric to a suitably defined (spatial) Schwarzschild manifold of dimension $n+1$ of positive mass.

\subsection{Asymptotically hyperbolic manifolds} \label{subsec-intro-AH}
Roughly speaking, a Riemannian manifold $(M^3,\g)$ is said to be \emph{asymptotically hyperbolic}, if, outside a compact set, it is diffeomorphic to $\Hb^{3} \setminus \ovr{B_1^{\Hb}(0)}$, and the metric $\g$ approaches the metric $\g_{H}=(1+r^{2})^{-1}dr^{2}+r^{2} g_{*}$ of the hyperbolic space $\Hb^{3}$ as $r\to\infty$. In this setting, the notion of total mass is more subtle and has been defined by X.~Wang~\cite{Wang} and by Chru\'sciel and Herzlich~\cite{CH}; here, it will be denoted it by $m_{\AH}(M^{3},\g)$. It vanishes on hyperbolic space, $m_{\AH}(\Hb^{3},\g_{H})=0$. The dominant energy condition for a time-symmetric initial data set can also be derived in this context --- after an appropriate rescaling of the so-called cosmological constant in the Einstein Equations --- and amounts to $R(\g) \geq -6$. The asymptotically hyperbolic Positive Mass Theorem (with rigidity $m_{\AH}(M^{3},\g)=0$ if and only if $(M^{3},\g)$ is isometric to hyperbolic space) holds for spin manifolds in all dimensions as was shown under restrictions on the asymptotics by X.~Wang~\cite{Wang} and without such restrictions by Chru\'sciel and Herzlich~\cite{CH}, and under some restrictions on the geometry at infinity as was shown by Andersson, Cai, and Galloway~\cite{ACG}, by Chru\'sciel, Galloway, Nguyen, and Paetz~\cite{CGNP}, and by Chru\'sciel and Delay~\cite{CD}. Sakovich~\cite{S} has announced a proof of the $3$-dimensional asymptotically hyperbolic Positive Mass Theorem in full generality (with extra asymptotic assumptions needed for asserting rigidity). The rigidity case has been settled under fully general asymptotic assumptions by Huang, Jang and Martin~\cite{HJM}.

The conjectured asymptotically hyperbolic Riemannian Penrose Inequality states that if $(M^{3},\g)$ is an asymptotically hyperbolic manifold with boundary $\pr M$ which consists of an outer-minimizing horizon, then
\begin{align} \label{eq-RPIAH-bound}
m_{\AH}(M^{3},\g) \geq \sqrt{\frac{|\pr M|_{\g}}{16\pi}}\( 1 + \frac{|\pr M|_{\g}}{4 \pi}\),
\end{align}
where $|\pr M|_{\g}$ represents the area of the boundary of $M$ with respect to $\g$ as before.

Equality holds if and only if $(M^{n+1},\g)$ is isometric to an AdS-Schwarzschild manifold of positive mass: The \emph{(spatial) AdS-Schwarzschild manifold of mass $m>0$} can be described as the manifold $(r_+,\infty) \times \bS^2$, with the Riemannian metric 
\begin{align*}
\g_{m,\AdS} =\( 1 + r^2 - \frac{2m}{r}\)^{-1} dr^2 + r^2 g_*,
\end{align*}
where $r_+$ is the largest positive zero of $1 + r^2 - \frac{2m}{r}$ and $g_*$ again denotes the standard round metric on $\bS^2$. This inequality has been proven in special cases by Dahl, Giquaud and Sakovich in~\cite{DGS}, complemented by the work of de Lima and Gir\~{a}o in \cite{dLG}; and by Ambrozio in \cite{AL}.

The \emph{hyperbolic Hawking mass} of a closed surface $(\S,g)$ inside $M^3$ with mean curvature $H$ is defined as
\begin{align*}
\m_\textnormal{H}^{\AH}(\S,g,H) \defeq \sqrt{\frac{|\S|_g}{16 \pi}}\( 1 - \frac{1}{16 \pi} \int_{\S} (H^2-4) \, d\s \).
\end{align*}

In~\cite{CCM}, mimicking the definition of Bartnik mass in the asymptotically flat case, McCormick and the authors defined a notion of  \emph{asymptotically hyperbolic Bartnik mass}. Given Bartnik data $(\S \cong \bS^2,g,H)$, we can define the set of admissible extensions~$\mathcal{A}$ as the set of asymptotically hyperbolic manifolds $(M^3,\g)$ with $R(\g) \geq -6$ such that their boundary $\pr M$ is outer-minimizing and isometric to $(\S,g)$ with mean curvature~$H$ (cf.~\cite{CCM}). The \emph{asymptotically hyperbolic Bartnik mass} is given by

\begin{align*}
\m_\textnormal{B}^{\AH}(\S,g,H) \defeq \inf\{ m_{\AH}(M^3,\g)\,\vert\,{(M^3,\g) \in \mathcal{A}}  \}.
\end{align*}

Clearly, if the asymptotically hyperbolic Riemannian Penrose conjecture is proven to be true, the hyperbolic Hawking mass would provide a lower bound for the hyperbolic Bartnik mass for an outer-minimizing connected minimal boundary.

\subsection{Initial data sets for the Einstein--Maxwell Equations} \label{subsec-intro-EM}

In the context of solutions to the Einstein--Maxwell Equations, i.e., when the Einstein equations are generalized to incorporate electromagnetic effects, and assuming that the magnetic field vanishes, a time-symmetric initial data set satisfying the dominant energy condition corresponds to a triple $(M^3,\g,E)$, where $(M^3,\g)$ is an asymptotically flat Riemannian manifold and $E$ is a divergence-free vector field on $(M^{3},\g)$, modeling the \emph{electric field}, satisfying
\begin{align*}
R(\g) \geq 2|E|_{\g}^2.
\end{align*}

An very important example is given by the \emph{(spatial) Reissner--Nordstr\"om manifold} of charge $Q$ and mass $m > |Q|$, is defined as the manifold $(r_+,\infty) \times \bS^2$, with the Riemannian metric given by
\begin{equation*}
\g_{m,Q} \defeq \( 1 - \frac{2m}{r} + \frac{Q^2}{r^2} \)^{-1}dr^2 + r^2 g_*, 
\end{equation*} 
where $r_+ \defeq m + \sqrt{m^2 - Q^2}$ and $g_*$ again denotes the standard round metric on~$\bS^2$. When $m > |Q|$ we call this a \emph{sub-extremal} spatial Reissner--Nordstr\"om manifold. The corresponding electric field is given by
\begin{equation*}
E_{m,Q} \definedas \frac{Q}{r^2}\sqrt{1- \frac{2m}{r} + \frac{Q^2}{r^2}}\, \pr_r.
\end{equation*}

The Riemannian Penrose Inequality in this context, called the \emph{Riemannian Penrose Inequality with (electric) charge}, was proven by Jang~\cite{Jang} assuming the existence of a smooth solution to the inverse mean curvature flow (IMCF); with the weak formulation of the IMCF by Huisken and Ilmanen~\cite{HI}, its general proof (for a connected outer-minimizing horizon) is now complete (see the works of Disconzi and Khuri~\cite{DK} and Mars \cite{Mars}). The corresponding rigidity statement was proven by Disconzi and Khuri in~\cite{DK}. The Riemannian Penrose inequality with (electric) charge can be stated as follows: Let $(M^3,\g,E)$ be a triple such that $(M^3,\g)$ is an asymptotically flat Riemannian manifold and $E$ is a divergence-free vector field on $M^{3}$ such that $R(\g) \geq 2|E|^2_{\g}$. Moreover, suppose that $\pr M$ consists of a connected, outer-minimizing minimal surface. Then
\begin{align} \label{eq-RPICH-bound}
m_{\ADM}(M^3,\g) \geq  \sqrt{\frac{|\pr M|_{\g}}{16 \pi}}\(1+ \frac{4\pi}{|\pr M|_{\g}}Q^2_{\infty}\),
\end{align}
where 
\begin{align*}
Q_{\infty}\definedas\lim_{r\to\infty}\frac{1}{4\pi}\int_{\bS^{2}_{r}} \g(E,N_{r})\,d\s_{r}
\end{align*}
denotes the \emph{total (electric) charge of $(M^{3},\g,E)$}, with $\bS^{2}_r\subset\R^{3}$ the coordinate sphere of radius $r$, $N_{r}$ the (Euclidean) outer unit normal to $\bS^{2}_{r}$, and $d\s_r$ the (Euclidean) area form on $\bS^{2}_r$. Equality holds if and only if $(M^{3},\g)$ is isometric to (an exterior region of) a sub-extremal Reissner--Nordstr\"om manifold of mass $m=m_{\ADM}(M^{3},\g)$ and charge $Q_{\infty}$, and $E$ corresponds to $E_{m,Q_{\infty}}$ (see \cite{Mars} and the review article by Dain and Gabach-Clement \cite{D-GC} for more detailed expositions of initial data sets for the Einstein--Maxwell Equations and related geometric inequalities).

Let $(\S,g)$ be a closed surface contained in $(M^3,\g,E)$. The \emph{total charge contained in~$\S$} is defined as $Q \defeq \frac{1}{4\pi}\int_{\S} \g(E,N)\,d\s$, where $N$ is the outward unit normal to $\S$ in $(M^{3},\g)$ and $d\s$ is the induced area form of $\S$. Let $H$ denote the mean curvature of $\S$ in $(M^3,\g)$, then the \emph{charged Hawking mass} of $(\S,g,H,Q)$ is given by
\begin{align*}
\m_\textnormal{H}^{\CH}(\S,g,H,Q) \defeq \sqrt{\frac{|\S|_{g}}{16 \pi}}\( 1 + \frac{4\pi Q^2}{|\S|_{g}} -  \frac{1}{16\pi} \int_{\S} H^2 \, d\s  \).
\end{align*}
\medskip

As in the other settings described above, it is possible to define a set of admissible extensions $\mathcal{A}$ for suitable charged Bartnik data, leading to an ad-hoc version of Bartnik mass in this scenario which was done by Alaee and the authors in \cite{ACC}. For given \emph{charged Bartnik data} $(\S \cong \bS^2,g,H,Q)$, where $(\S,g,H)$ is as above and $Q \in \R$, we say that a triple $(M^{3},\g,E)$ is an admissible extension of $(\S,g,H,Q)$ if $(M^{3},\g)$ is an asymptotically flat Riemannian manifold with outer-minimizing boundary $\pr M$ isometric to $(\S,g)$ with mean curvature $H$ in $(M^{3},\g)$ and $E$ is a smooth vector field on $M^{3}$, interpreted as an electric field, such that $R(\g) \geq |E|^2_{\g}$ and $\dv_{\g} E =0$, and such that $Q$ is the total charge in $\pr M$. Then, we define the \emph{charged Bartnik mass} as
\begin{align*}
\m_\textnormal{B}^{\CH}(\S,g,H,Q) \defeq \inf\{ m_{\ADM}(M^3,\g)\,\vert\,{(M^3,\g) \in \mathcal{A}}  \}.
\end{align*}

\section{The Mantoulidis--Schoen Construction} \label{sec-MS}

\begin{figure}[ht] 
\includegraphics{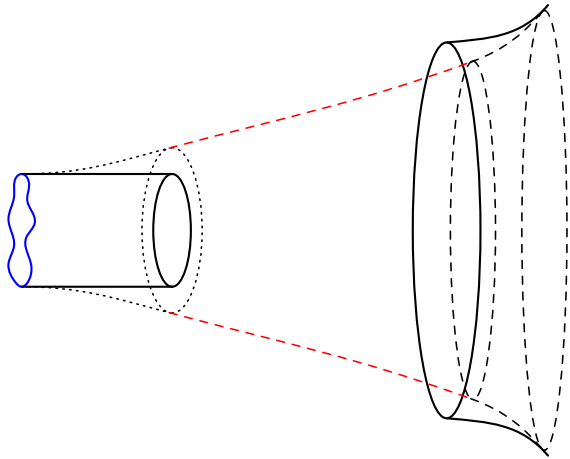}
\caption{The Mantoulidis--Schoen construction~\cite{MS}. The {\color{blue}blue} surface represents the given surface $(\S,g)$, dotted $3$-manifold represents the collar extension, the black dashed $3$-manifold represents the defomed Schwarzschild manifold of mass $m$, the {\color{red} red} dashed $3$-manifold represents the positive scalar bridge.} \label{fig-MS}
\end{figure}

Recall that, as discussed in Section~\ref{sec-intro}, the Bartnik mass is remarkably difficult to compute. In~\cite{MS}, Mantoulidis and Schoen computed the Bartnik mass for minimal Bartnik data --- $(\S \cong \bS^2, g, H)$ is \emph{minimal Bartnik data} if $H \equiv 0$. To do this, given $(\S \cong \bS^2, g, H \equiv 0)$ satisfying $\lambda_1 \defeq \lambda_1(-\Lap_g + K(g))>0$, where $K(g)$ is the Gaussian curvature and $\lambda_1$ denotes the first eigenvalue of $-\Lap_g + K(g)$, they constructed asymptotically flat extensions with non-negative scalar curvature, such that the ADM mass of the extensions can be made arbitrarily close to the optimal value of~\eqref{eq-RPI-bound}. The condition on $\lambda_1$ arises naturally in the theory of minimal surfaces, and its connection to Mathematical Relativity comes from the fact that outer-minimizing horizons in an asymptotically flat initial data set with non-negative scalar curvature are, in fact, stable minimal spheres. The construction can be described by the following four steps:
\begin{enumerate}
\item Using the uniformization theorem, write $g=e^{2w} g_*$, where $w$ is a smooth function on $\S$. Connect $g$ to $g_*$ by a smooth path of metrics given by $g(t)=e^{2(1-t)w} g_*$, $t\in[0,1]$.
\item Modify $\{ g(t) \}_{0 \leq t \leq 1}$ by a $t$-dependent family of diffeomorphisms on $\S$ in a suitable way, and use it to construct a \emph{collar extension} of $(\S,g)$ with positive scalar curvature, that is, a Riemannian 3-manifold with topology $[0,1] \times \S$, such that $\S_0 = \{0\} \times \S$ is isometric to $(\S,g)$ and $\S_1=\{ 1 \} \times \S$ is round. This collar extension is constructed so that  $\S_0$ is minimal and outer-minimizing. Additionally, it is rotationally symmetric in a region near $\S_1$, and $|\S_1|_{g(1)}$ can be made arbitrarily close to $|\S_0|_{g}$. See dotted black part in Figure~\ref{fig-MS}.  \newpage
\item For any $m > \sqrt{\frac{|\S|_{g}}{16 \pi}}$, consider a Schwarzschild manifold of mass $m$ and deform it in a rotationally symmetric way in a small region near the horizon so that it has positive (instead of vanishing) scalar curvature in this region. See dashed black part in Figure~\ref{fig-MS}.
\item Smoothly glue these two manifolds via a positive scalar curvature bridge. See dashed red part in Figure~\ref{fig-MS}. 
\end{enumerate}

The resulting manifold is a smooth, asymptotically flat manifold with non-negative scalar curvature, such that its boundary is isometric to $(\S,g)$, minimal and outer-minimizing, and hence, an admissible extension of $(\S,g,H\equiv 0)$. Since this manifold is exactly Schwarzschild of mass $m$ outside a compact set, its ADM mass is equal to $m$, which can be made arbitrarily close to the optimal value in~\eqref{eq-RPI-bound}, proving that 
\begin{align*}
m_B(\S,g,H \equiv 0) = \sqrt{\frac{|\S|_g}{16 \pi}} = \m_\textnormal{H}(\S,g,H\equiv 0).
\end{align*}

Besides providing an exact value for the Bartnik mass of suitable minimal Bartnik data, the Mantoulidis--Schoen construction can be thought as evidence for instability of the Riemannian Penrose Inequality (Theorem~\ref{thm-RPI}) in the sense that, given $(\S \cong \bS^2,g)$, the ADM mass of the extension can be made arbitrarily close to the optimal value of \eqref{eq-RPI-bound}, while the geometry of the horizon could be far away from being round. For a more detailed summary of the results of Mantoulidis and Schoen~\cite{MS} and interesting consequences thereof, see~\cite{MS-CQGplus} by the same authors.

From this perspective and from considering the Riemannian Penrose Inequality in higher dimensions \eqref{eq-RPI-hd}, it is natural to ask whether a construction similar to this one could be performed in higher dimensions. Miao and the first named author answered this question positively in~\cite{CM}, by extending the arguments developed in~\cite{MS} to higher dimensions. More precisely, they showed that given $(\S \cong \bS^{n},g,H \equiv 0)$ $(n \geq 3)$ with $R(g) > 0$ for $n=3$ and for $n \geq 4$ with the additional assumption that $(\S,g)$ can be isometrically embedded into $(\R^{n+1},\delta)$ as a star-shaped surface, then, for any $m>0$ such that
\begin{align*}
m > \frac{1}{2} \( \frac{|\S|_g}{\w_n} \)^{\frac{n-1}{n}},
\end{align*}
where $\w_n$ denotes the volume of $\bS^n$ and $|\S|_g$ the volume of $(\S,g)$, there exists an asymptotically flat extension with non-negative scalar curvature and ADM mass equal to $m$, such that its boundary is isometric to $(\S \cong \bS^{n},g)$, minimal and outer-minimizing. A key step in this construction is the existence of paths connecting the given metric $g$ to a round metric (corresponding to Step (1) above). For $n=3$, this is possible since the space of metrics with positive scalar curvature on $\bS^3$ is path connected, which was proved by Marques using Ricci flow with surgeries in~\cite{MF}. For $n \geq 4$, extrinsic geometric flows developed by Gerhardt~\cite{GC} and Urbas~\cite{UJ} can be used.

\begin{figure}[h] 
\includegraphics{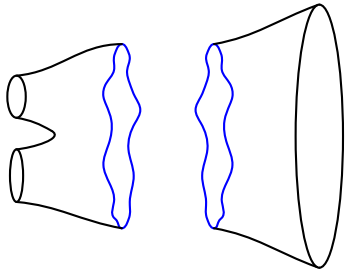}
\caption{The Miao--Xie construction~\cite{MX}. The left part depicts the manifold $\Omega$, in which the CMC boundary component $\S$ is depicted in {\color{blue} blue}. The right part depicts the collar extension of~$\S$.} \label{fig-MX} 
\end{figure}

Motivated by the quasi-local mass problem in General Relativity, in~\cite{MX}, Miao and Xie investigated the effect of non-negative scalar curvature in compact bodies inside a Riemannian manifold. More precisely, they considered compact 3-dimensional Riemannian manifolds $\Omega$ with boundary $\pr \Omega= \S_h \cup \S$, where $\S \cong \bS^2$ has positive constant mean curvature $H$, and $\S_h$ is a minimal surface (possibly disconnected), and such that there are no other minimal surfaces in $\Omega$. Inspired by~\cite{MS}, they constructed collar extensions of $(\S,g)$ where $g$ is a non-round metric induced on $\S$ by $\Omega$, as in Step (2) above, but with the property that the mean curvature of the boundary $\S$ is equal to $H$. Then, they attach this collar extension to $\Omega$ along $\S$ and apply a result by Miao~\cite{Miao-locRPI} to obtain
\begin{align*}
\m_\textnormal{H}(\S_1) \geq \sqrt{\frac{|\S_h|}{16 \pi}}.
\end{align*}   

Assuming a smallness condition on $H$ (only in terms of $(\S,g)$), they were able to estimate $\m_\textnormal{H}(\S_1)$ in terms of $\m_\textnormal{H}(\S)$ and quantities that depend only on $(\S,g)$, which translates to relations between $\S_h$ and $\S$, and information about $\S$ itself when $\S_h = \emptyset$. This construction and its consequences provide insights about the interaction of the scalar curvature and  geometry of the vicinity of a horizon in a Riemannian $3$-manifold. We refer to the interested reader to~\cite{MX} for the precise statements and a more detailed description. The collar extensions used by Miao and Xie in~\cite{MX}, which will be discussed in Section~\ref{sec-Bartnik-est}, have the remarkable property that the Hawking mass of $\S_t=\{ t \} \times \S$ is controlled along the collar. 

\begin{figure}[ht] 
\includegraphics{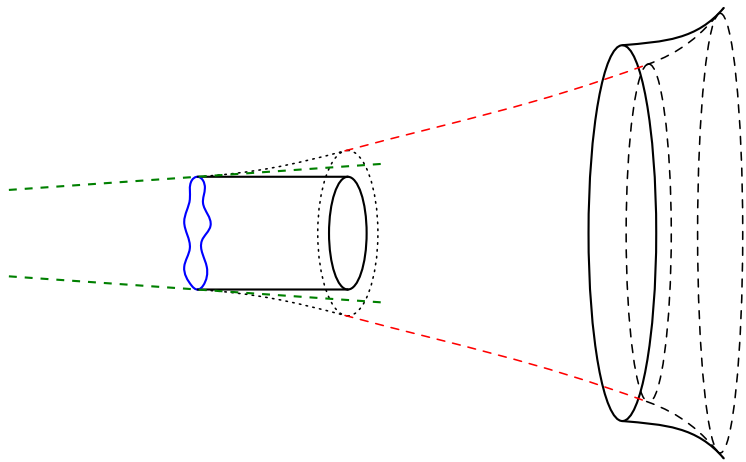}
\caption{The modification of the Mantoulidis--Schoen construction to allow a CMC boundary in~\cite{CCMM}. The given surface $(\S,g)$ is represented in {\color{blue} blue} and the given constant mean curvature is represented by the dashed {\color{green!50!black}green} line. As before, the collar extension is represented by the dotted part and the deformed Schwarzschild manifold by the dashed black line. The positive scalar bridge is represented by the {\color{red}red} dashed part.} \label{fig-CCMM}
\end{figure}

In~\cite{CCMM}, McCormick, Miao, and the authors developed analytic tools to smoothly glue compact, rotationally symmetric manifolds with boundary $\Omega_{\textnormal{rot}}$ to (an exterior region of) a Schwarzschild manifold, inspired by the tools developed in~\cite{MS}. The main difference is that in~\cite{CCMM}, the set of Schwarzschild manifolds for which this gluing technique works depends on the Hawking mass of the ``outer'' boundary of~$\Omega_{\textnormal{rot}}$, instead of on its area. Since the collar extensions constructed by Miao and Xie in~\cite{MX} give a very precise control on the Hawking mass growth along the collar, and since by Step (2) above they are rotationally symmetric close to the ``outer'' boundary, combined with the gluing construction in~\cite{CCMM}, they produce asymptotically flat extensions of a given surface $(\S,g)$ such that its mean curvature in the extension is equal to a  prescribed, sufficiently small positive constant $H$. Furthermore, the resulting manifolds are admissible extensions in the sense of Section~\ref{sec-intro}, so this result provided new estimates for \emph{constant mean curvature (CMC)} Bartnik data (see Figure~\ref{fig-CCMM}).

After the success of the Mantoulidis--Schoen construction and its modifications~\cite{MX,CCMM} in computing the Bartnik mass for minimal Bartnik data and in estimating the Bartnik mass of CMC Bartnik data in the context of asymptotically flat manifolds, McCormick and the authors investigated in~\cite{CCM} if similar constructions could be carried out in the context of asymptotically hyperbolic manifolds. After formulating an analog to the Bartnik mass in the asymptotically hyperbolic case, see Section~\ref{subsec-intro-AH}, similar results for both the minimal and the CMC case are obtained, following the techniques in~\cite{MS,MX,CCMM}.

%Very recently, the techniques of~\cite{CCM} were combined by Gehring in~\cite{G} with those developed in~\cite{CM}, to handle the minimal case in the asymptotically hyperbolic setting in higher dimensions $n\geq3$.

Later, Miao, Y.~Wang, and Xie~\cite{MWX} refined the construction of the collar extensions by Miao and Xie~\cite{MX}, removing the smallness condition on $H$. In particular, this refinement allowed them to find conditions on the intrinsic metric of the boundary $\S$ of a domain $\Omega$ in a Riemannian manifold $(M^3,\g)$ with non-negative scalar curvature that guarantees the positivity of the Hawking mass of $\S$, previously known for stable CMC surfaces by work of Christodoulou and Yau~\cite{CY}. In view of Bartnik mass estimates for CMC Bartnik data $(\S \cong \bS^2,g,H)$, the refined collar construction allowed Miao, Y.~Wang, and Xie~\cite{MWX} to estimate the Bartnik mass without the smallness assumption on $H$ needed in~\cite{CCMM}, instead simply requiring the Hawking mass of the given Bartnik data to be non-negative. They also improved the collar construction and obtained related Bartnik mass estimates in the asymptotically hyperbolic case.

Very recently, Alaee and the authors formulated an ad-hoc definition of a charged Bartnik mass in the context of the Einstein--Maxwell Equations~\cite{ACC}, see Section~\ref{subsec-intro-EM}, and computed its value for charged minimal Bartnik data. This construction and computation were applied and generalized to higher dimensions by Pe\~nuela D\'iaz in \cite{P}, utilizing techniques and results of~\cite{CM}, and using the Ricci flow existence and convergence results for compact Riemannian manifolds with pointwise $\tfrac{1}{4}$-pinched sectional curvature developed by Brendle and Schoen in~\cite{BS,BS2}.

Also very recently, the techniques of~\cite{CCM} were combined by Gehring in~\cite{G} with those developed in~\cite{CM}, to handle the minimal case in the asymptotically hyperbolic setting in higher dimensions $n\geq3$.

The Mantoulidis--Schoen construction, its parts and modifications, have been used  for various purposes besides constructing admissible extensions for the Bartnik mass. In this survey, we will only focus on the applications mentioned above. However, in Section~\ref{sec-conclusions} we provide a list of references to works where aspects of this construction were used.

\section{Collar Extensions, Gluing Procedures, and Bartnik Mass Estimates} \label{section-gluingmethods}

We will now describe the Mantoulidis--Schoen construction and its consequent extensions to more general settings in some detail. For the sake of presentation, we provide the descriptions and computations in the most general form known to us, which is a result of the contributions made in the articles mentioned in Section~\ref{sec-MS}. To see the precise statements, we invite the reader to consult those works individually.

This section is divided into three parts. The first one is devoted to describe the collar extensions in general, the second one to explain the gluing techniques in different contexts, and the last one to set the previous two parts in the context of Bartnik mass estimates. While the first two are relevant in higher dimensions $n \geq 2$ and may well have or are known to have applications in other settings, Bartnik mass estimates will only be addressed for $n=2$ because of the predominantly physically motivated interest in quasi-local notions of mass.

\subsection{Collar extensions}

Let $(\S \cong \bS^{n},g)$ be an $n$-dimensional Riemannian manifold. A \emph{collar extension} of $\S$ is an $(n+1)$-dimensional Riemannian manifold $(M,\g)$, where $M=[0,1] \times \S$, with a Riemannian metric of the form
\begin{equation} \label{eq-metric-c}
\g \defeq v(t,\cdot)^2 dt^2 + h(t),
\end{equation}
where $v$ is a positive smooth function on $[0,1] \times \S$  and $\{ h(t) \}_{0 \leq t \leq 1}$ is a smooth family of metrics on $\S$, satisfying $h(0)=g$. Here, we chose the interval $[0,1]$ for convenience. Notice that the condition $h(0)=g$ prescribes the metric on what we call the \emph{inner boundary}, making it isometric to $(\S,g)$.

Denote $\frac{d}{dt}$ by $'$. Our convention for the mean curvature is that spheres in Euclidean space have positive mean curvature with respect to the outward unit normal. The \emph{outward} unit normal to $\S_t$ in a collar extension is the one pointing in the direction of $\pr_t$. Hence, a direct computation shows that the mean curvature for the level set $\S_t=\{ t \} \times \S$, is given by
\begin{equation} \label{eq-mean-c}
H \defeq H(t,\cdot) = \frac{1}{2 v(t,\cdot)} \tr_{h(t)} h'(t).
\end{equation}
This is one of the very useful relations between the family of metrics $\{ h(t) \}_{0 \leq t \leq 1}$ and the geometry of the collar extension $(M,\g)$. For example, if we would like $(M,\g)$ to have a minimal inner boundary, then we would need to arrange $\tr_{h(t)} h'(t) \vert_{t=0}=0$.

We now proceed to compute the scalar curvature $R(\g)$ of the collar extension~\eqref{eq-metric-c}. Denoting the scalar curvature and the second fundamental form of $\S_t$ in $M$ by $R(h)$  and $\II\defeq \II(t,\cdot)$, respectively, and the outward unit normal to $\S_t$ by $N = N(t,\cdot)$, using the Gauss equation we know that 

\begin{align*}
R(\g) &= 2\Ric_{\g}(N,N) + R(h) - H^2 + |\II|^2, \\
&= 2\Ric_{\g}(N,N) + R(h) - \frac{1}{4v^2} \( \tr_h h'    \)^2 + \frac{1}{4 v^2} |h' |^2_h.
\end{align*}
Direct computations give
\begin{align*}
\Ric_{\g}(N,N) = - 2\frac{\Lap_h v}{v}   -  \frac{\tr_{h}{h''}}{v^{2}} + \frac{\pr_t v\tr_h h'}{v^{3}} + \frac{1}{2v^{2}} |h'|^2_h.
\end{align*}
Therefore, we have
\begin{align} 
\begin{split} \label{eq-scalar-collar-gen}
R(\g) &=  \frac{2}{v}\(-\Lap_h v + \frac{R(h)}{2}v\)   \\ 
&\qquad + \frac{1}{v^2}\(- \tr_{h}{h''} + \frac{\pr_t v}{v}\tr_h h'   
 - \frac{1}{4} \( \tr_h h'    \)^2 + \frac{3}{4 } |h' |^2_h  \).
\end{split}
\end{align}

We are interested in the case where $h(t) = F(t)^2 g(t)$, where $F(t) \geq  1$ is the \emph{radial profile} of the collar extension, which is a smooth function on $[0,1]$ with $F(0)=1$, and $\{ g(t) \}_{0 \leq t \leq 1}$ is a smooth family of metrics satisfying\label{i-iv}
\begin{enumerate}[(i)]
\itemsep0.25em
\item $g(0)=g$ and $g(1)$ is round,
\item $g'(t) \equiv 0$ for $t \in [\theta,1]$ for some $0 < \theta < 1$, and
\item $\tr_{g(t)} g'(t) \equiv 0$ for all $t \in [0,1]$.
\end{enumerate}
Depending on the situation, the smooth family of metrics $\{ g(t) \}_{0 \leq t \leq 1}$ will be required to satisfy given special curvature conditions. If we denote the set of all metrics satisfying the respective given special curvature conditions by $\M$, the latter can then be expressed as
\begin{enumerate}[(i)] 
\setcounter{enumi}{3}
\item $g(t) \in \M$  for all $t \in [0,1]$.
\end{enumerate}

%\newpage

\begin{remark}
Property (i) prescribes the metric at the boundaries of the collar extension, while (ii) implies that $(M,\g)$ is rotationally symmetric for $t \in [\theta,1]$. The property (iii) can be seen to be equivalent to requiring that the volume form $dV_{g(t)}$ is preserved along the path $\{ g(t) \}_{0 \leq t \leq 1}$.
\end{remark}
\newpage
\begin{remark}\label{rem:ab}
It is possible to associate two scaling-invariant quantities to the path constructed above. These quantities can be regarded as a measurement of how round $(\S,g)$ is, and they are given by 
\begin{align}\label{eq-alpha}
\a\defeq \frac{1}{4}\max_{[0,1] \times \S} | g'(t)|^2_{g(t)}, 
\end{align}
and
\begin{align}\label{eq-beta} 
\b \defeq \min_{[0,1] \times \S} r_o^n \frac{R(g(t))}{n(n-1)},
\end{align}
where $r_o \defeq \left(\frac{|\S|_{g}}{\w_{n}}\right)^{\frac{1}{n}}$ is the \emph{volume radius} of $(\S,g)$. Note that by (iii) this implies that $|\S|_{g(t)} \equiv \w_{n} r_o^n$ for all $t \in [0,1]$. Note that if $g$ is round then a constant path can be chosen, and for this path, $\a=0$ and $\b=1$. The quantities $\a$ and $\b$ were defined and studied in detail by Miao and Xie~\cite{MX} in dimension $n=2$.
\end{remark}

The existence of a path satisfying (i)-(iii) follows after we know that a volume-preserving smooth family of metrics $\{ \widetilde{g}(t) \}_{0 \leq t \leq 1} \subset \M$ satisfying (i) exists: Indeed, by reparametrizing the family $\{ g(t) \}_{0 \leq t \leq 1}$ through a composition with a bump function, one can achieve (ii) (which of course preserves (i) and (iv)). Subsequently, one can find a 1-parameter family of diffeomorphisms $\{ \phi_t \}$ on $\S$ such that $g(t) \defeq \phi_t^*(\widetilde{g}(t))$ satisfies (iii), by solving a 1-parameter family of Poisson problems on $(\S,\widetilde{g}(t))$ (in~\cite{MS}, the method was inspired by what is called ``Moser's trick" in symplectic geometry). Note that this does not affect the validity of (i), (ii), and (iv). See~\cite[Lemma 1.2]{MS} for details when $n=2$ and~\cite[Lemma 4.1]{CM} for the extension to arbitrary dimensions.

With this choice of $\{ h(t) \}_{0 \leq t \leq 1}$ (i.e., $h(t)=F(t)^2g(t)$ and $\{ g(t) \}_{0 \leq t \leq 1}$ satisfying (i)--(iv) above), we have
\begin{align*}
\tr_{h} h' &= 2n\frac{F'}{F} + \tr_{g}{g'}= 2n\frac{F'}{F} , \\
|h'|_{h}^2 &= 4n\frac{(F')^2}{F^2} +4\frac{F'}{F}\tr_{g}g' + |g'|^2_{g}= 4n\frac{(F')^2}{F^2} + |g'|^2_{g}, \textnormal{and}\\
\tr_{h}h'' &= 2n\left[\frac{(F')^2}{F^2} + \frac{F''}{F}\right]  +
 4 \frac{F'}{F}\tr_{g}g' + \tr_{g}{g''} \\
 &=2n\left[\frac{(F')^2}{F^2} + \frac{F''}{F}\right]  + |g'|^2_{g}, 
\end{align*}
since $\tr_{g}g' \equiv 0$ by (iii). Using \eqref{eq-scalar-collar-gen}, the scalar curvature of the collar extension is given by
\begin{align} 
\begin{split}\label{eq-scal-collar}
R(\g) &=  \frac{2}{F^2v}\(-\Lap_{g} v + \frac{R(g)}{2}v\)  \\
& \qquad+ \frac{1}{v^2}\(2n\frac{F'}{F}\frac{\pr_t v}{v}  - \frac{n}{F^2}\left[(n-1)F'^2 + 2FF''\right] -\frac{1}{4 } |g'|^2_{g}
 \).
\end{split}
\end{align}
Furthermore, from \eqref{eq-mean-c}, we find that the mean curvature of $\S_t$ is given by
\begin{align}\label{eq-meanc-leaf}
H(t)= \frac{n}{ v(t,\cdot)} \frac{F'(t)}{F(t)}.
\end{align}

For the remainder of this survey, we assume that $h(t)=F(t)^2g(t)$ in \eqref{eq-metric-c} as described above.

This technique for constructing collar extensions of given data $(\S \cong \bS^n,g)$ corresponds to Step (2) in the Mantoulidis--Schoen construction. In what follows, we will explain how to glue the collar manifold to an exterior region of a rotationally symmetric model manifold . When this manifold is a $3$-dimensional Schwarzschild manifold of positive mass, this corresponds to Steps (3) and (4) in the Mantoulidis--Schoen construction.

\subsection{Smooth gluing procedures} \label{subsec-gluing}

A collar extension as in \eqref{eq-metric-c} with $h(t)=F(t)^2g(t)$, where the smooth path of metric $\{ g(t) \}_{0 \leq t \leq 1}$ satisfies (i)-(iv), has the special feature that  it is rotationally symmetric for $t \in [\theta,1]$, provided that $v(t,\cdot)$ is constant on $\S_t$ for $t \in [\theta,1]$ (which will be the case for all choices of $v$ considered here). Since Schwarzschild manifolds and all other model manifolds considered here are rotationally symmetric, Step (4) in the Mantoulidis--Schoen construction and its various generalizations can thus be phrased as a problem concerning how to smoothly glue two rotationally symmetric manifolds.

\begin{figure}[h] 
\includegraphics{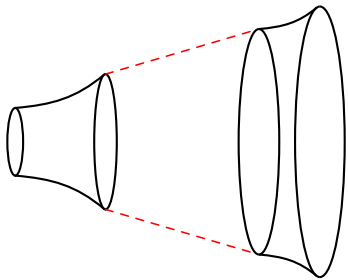}
\caption{The black solid part represent two rotationally symmetric Riemannian manifolds. The {\color{red}red} dashed part represents a piece of a cone attached between them. The joining points will not be smooth in general.} \label{fig-ConG}
\end{figure}

More precisely, consider two rotationally symmetric Riemannian manifolds $(M_1,\g_1)=([a_1,b_1] \times \bS^{n},ds^2 + f_1(s)^2 g_*)$ and $(M_2,\g_2)=([a_2,b_2] \times \bS^{n},ds^2 + f_2(s)^2 g_*)$, both with non-negative scalar curvature. Is there a way to smoothly glue these two manifolds in such a way that the non-negativity of the scalar curvature is preserved? Without losing generality, we may assume that $f_1(b_1) \leq f_2(a_2)$ and $b_1 < a_2$. Then, a natural approach will be to attach them together via a piece of a cone (see Figure~\ref{fig-ConG}). Note that, in general, the resulting manifold would of course not be smooth. To smooth out this manifold, one could perform a mollification, however, this could disrupt the non-negativity of the scalar curvature (say, if the scalar curvature at $s=a_1$ or $s=a_2$ is equal to 0, as it is for a Schwarzschild manifold).

To make some room for the mollification to work, we impose the condition that $M_1$ and $M_2$ have strictly positive scalar curvature. In~\cite{MS}, Mantoulidis and Schoen proved an analytic tool to smoothly glue two rotationally symmetric manifolds with positive scalar curvature for $n=2$ (under some extra conditions that we will state in a moment). This tool was later generalized to other settings, see below.\\

Observe that a metric of the form $\g=v_0^2 dt^2 + F_0(t)^2 g_*$ on $[a_0,b_0]\times \bS^{n}$, where $v_0$ is a positive constant and $F_0$ is a smooth positive function on  $[a_0,b_0]$ (which will be the case for the collar extensions described above on $[\theta,1]$, see also Section~\ref{sec-Bartnik-est}) can be written in the form $\g=ds^2+f(s)^2g_*$ after the change of variables $s=v_0t$, with $f(s)\definedas F(\frac{s}{v_{0}})$. Hence, the following lemma can be applied to the discussed collar extensions.

\begin{lemma}[Gluing Lemma~{\cite[Lemma 2.2]{MS}}] \label{lemma-gluing-MS}
Let  $f_i\colon [a_i,b_i] \to \R$ be two smooth functions ($\,i=1,2$). Suppose that 
\begin{enumerate}[(a)]
\itemsep0.25em
\item $f_i,f_i',f_i''>0$ for $i=1,2$,
\item the scalar curvature of the metric $\g_i = ds^2 + f_i(s)^2 g_*$ on $[a_i,b_i] \times \bS^2$ is positive for $i=1,2$, and
\item $f_1(b_1) < f_2(a_2)$ and $f_1'(b_1)=f_2'(a_2)$,
\end{enumerate}
then, after translating the intervals so that  $(a_2-b_1)f_1'(b_1)=f_2(a_2)-f_1(b_1)$, there exists a smooth function $f\colon [a_1,b_2]\to\R$ such that
\begin{enumerate}[(A)]
\itemsep0.25em
\item $f,f'>0$ on $[a_1,b_2]$,
\item $f \equiv f_1$ on $\left[ a_1, \frac{a_1+b_1}{2}\right]$, and $f\equiv f_2$ on $\left[ \frac{a_2+b_2}{2},b_2 \right]$, and
\item the metric $\g_f \defeq ds^2 + f(s)^2 g_*$ has positive scalar curvature on $[a,b] \times \bS^2$. 
\end{enumerate}
\end{lemma}

The Schwarzschild manifold~\eqref{eq-Schwarzschild-metric} can smoothly be extended to the boundary in such a way that we can write it as the manifold $[0,\infty) \times \bS^2$ with metric given by
\begin{align} \label{eq-Schwarzschild-metric2}
\g_m = ds^2 + u_m(s)^2 g_*,
\end{align}
where $u_m$ is a smooth function such that
\begin{enumerate}[(a)]
\itemsep0.25em
\item $u_m(0)=2m$,
\item $u_m'(s)=\sqrt{1 - \frac{2m}{u_m(s)}}$, and
\item $u_m''(s)=\frac{m}{u_m(s)^2}$.
\end{enumerate}

However, recall that the scalar curvature of a Schwarzschild manifold is zero, thus one can not directly apply Lemma~\ref{lemma-gluing-MS}. It turns out that it is possible to perform a deformation in a small region in a Schwarzschild manifold in order to push up the scalar curvature a bit~\cite[Lemma 2.3]{MS}. After this deformation, Mantoulidis and Schoen adjust the parameters of the collar extension and the deformation so that the conditions in Lemma~\ref{lemma-gluing-MS} are satisfied. Consequently, by considering collar extensions for which the area of $\S_1=\{1\} \times \S$ is not much bigger than the area of $\S_0=\{ 0 \} \times \S$, they obtain the desired manifold~\cite[Proof of Theorem 2.1]{MS}. This procedure was extended to higher dimensions by Miao and the first named author in~\cite{CM}, with appropriate generalizations of all relevant notions, see also Section~\ref{subsec-intro-higherdim}.

Following the general idea in~\cite{MS}, in the context of CMC Bartnik data and $n=2$, McCormick, Miao, and the authors developed a (smooth) gluing tool by relaxing the conditions in Lemma~\ref{lemma-gluing-MS} and expressing the conditions needed for the gluing in terms of the Hawking mass of the surface $\S_1$ in the collar extension. This defines a set of Schwarzschild manifolds that can be glued to a given rotationally symmetric manifold. For more details, in particular for an explanation why the size of the Hawking mass of $\S_1$ becomes relevant, see Section~\ref{subsec-CMC}. We state this smooth gluing tool in the following proposition.

\begin{prop}[CMC-Gluing to a Schwarzschild manifold~{\cite[Prop.~2.1]{CCMM}}] \label{prop-gluing-CCMM}
Consider a metric $\g_f\defeq ds^2 + f(s)^2 g_*$ on $[a,b] \times \bS^2$, where $f$ is a smooth, positive, and increasing function on $[a,b]$. Suppose that
\begin{enumerate}[(a)]
\itemsep0.25em
\item $\g_f$ has positive scalar curvature,
\item $\S_b = \{ b \} \times \bS^2$ has positive mean curvature, and
\item $\m_\textnormal{H}(\S_b) \geq 0$.
\end{enumerate}
Then, for any $m > \m_\textnormal{H}(\S_b)$, the manifold $([a,b] \times \bS^2,\g_f)$ can be smoothly glued to (an exterior region of) a Schwarzschild manifold of mass $m$, such that the resulting asymptotically flat manifold is rotationally symmetric, has non-negative scalar curvature and its coordinate spheres have positive constant mean curvature.
\end{prop}

Even though Proposition~\ref{prop-gluing-CCMM} is essentially a combination of Steps (3) and (4)  of the original construction, it has the advantage of avoiding a parameter matching procedure to perform the gluing, hence, it is suitable for applications in other situations, since it gives quantitative conditions in terms of the Hawking mass of $\S_b= \{ b \} \times \bS^2$ which allows to perform the gluing.

Later, McCormick and the authors~\cite{CCM} developed a gluing tool similar to Proposition~\ref{prop-gluing-CCMM} to glue rotationally symmetric manifolds to AdS-Schwarzschild manifolds, also for $n=2$, (cf. Section~\ref{subsec-intro-AH}), which will be presented next.

\begin{prop}[Gluing to an AdS-Schwarzschild manifold~{\cite[Prop.~3.3]{CCM}}]\label{prop-gluing-CCM}
Consider a metric $\g_f\defeq ds^2 + f(s)^2 g_*$ on $[a,b] \times \bS^2$, where $f$ is a smooth, positive, and increasing function on $[a,b]$. Suppose that
\begin{enumerate}[(a)]
\itemsep0.25em
\item $\g_g$ has scalar curvature $R(\g_f) > -6$,
\item $\S_b$ has positive mean curvature, and
\item $\m_\textnormal{H}^{\AH}(\S_b) \geq -f(b)^3$.
\end{enumerate}
Then, for any $m > \m_\textnormal{H}^{\AH}(\S_b)$, the manifold $([a,b] \times \bS^2,\g_f)$ can be smoothly glued to (an exterior region of) an AdS-Schwarzschild manifold of mass $m$, such that the resulting asymptotically hyperbolic manifold has scalar curvature bigger than or equal to $-6$ and its coordinate spheres have positive constant mean curvature.
\end{prop}

Following~\cite{CM,CCM}, this proposition has very recently been extended to higher dimensions by Gehring in~\cite{G}, with appropriate generalizations of all relevant notions, see also Sections~\ref{subsec-intro-higherdim} and \ref{subsec-intro-AH}.\\

In the context of asymptotically flat manifolds $(M^{3},\g)$ together with a vector field $E$ on $M^{3}$ playing the role of an electric field, that is, to systems $(M^{3},\g,E)$ corresponding to time-symmetric initial data sets for the Einstein--Maxwell Equations (see Section~\ref{subsec-intro-EM}), Alaee and the authors~\cite{ACC} developed the next proposition which provides a gluing tool to smoothly glue rotationally symmetric Riemannian manifolds with an electric field to a Reissner--Nordstr\"om manifold of positive mass and sufficiently small charge $Q$. The main difference between the following result and the previous ones relies on the extra work needed to also obtain a smooth extension of the electric field $E$ satisfying the desired properties.

\begin{prop}[Gluing to a Reissner--Nordstr\"om manifold~{\cite[Prop.~4.1]{ACC}}] \label{prop-gluing-ACC}
Consider the metric $\g_f\defeq ds^2 + f(s)^2 g_*$ on $[a,b] \times \bS^2$, where $f$ is a smooth, positive, and increasing function on $[a,b]$. Suppose that
\begin{enumerate}[(a)]
\itemsep0.25em
\item the scalar curvature of $\g_f$ satisfies
\begin{align*}
R(\g_f) > 2|E|^2_{\g_f},
\end{align*}
where $E\definedas Qf^{-2} \pr_s$ for some constant $Q \in \R$,
\item $f(b) > |Q|$,
\item $\S_b$ has positive mean curvature, and
\item $\m_\textnormal{H}^{\CH}(\S_b) \geq |Q|$.
\end{enumerate}
Then, for any $m > \m_\textnormal{H}^{\CH}(\S_b)$, the manifold $([a,b] \times \bS^2,\g_f)$ can be smoothly glued to (an exterior region of) a sub-extremal  Reissner--Nordstr\"om manifold of mass $m$ and charge $Q$. Moreover, the vector field $E$ can be smoothly extended to a  divergence free vector field $\widetilde{E}$, which eventually coincides with the electric field of the Reissner--Nordstr\"om manifold. In addition, the resulting asymptotically flat manifold has scalar curvature bigger than or equal to $2|\widetilde{E}|^2$ and its coordinate spheres have positive  constant mean curvature.
\end{prop}

Following~\cite{CM,ACC}, this proposition has very recently been extended to higher dimensions by Pe\~nuela D\'iaz in~\cite{P}, with appropriate generalizations of all relevant notions, see also Sections~\ref{subsec-intro-higherdim} and \ref{subsec-intro-EM}.\\

It can be shown using maximum principle methods from minimal and CMC surface theory that Riemannian manifolds $(M^{n},\g)$ foliated by positive constant mean curvature spheres automatically have outer-minimizing inner boundary. This allows us to use the above Lemmas and Propositions towards constructing admissible extensions for Bartnik data as introduced in Section~\ref{sec-intro}.

\subsection{Obtaining Bartnik mass estimates}

Suppose that $(\S \cong \bS^2,g,H)$ are given $2$-dimensional Bartnik data. It follows from the definition of Bartnik mass that in order to obtain an upper bound for $\m_\textnormal{B}(\S \cong \bS^2,g,H)$, it suffices to construct any admissible extension of the data --- in general, constructing admissible extensions is a very difficult problem. The Mantoulidis--Schoen construction has proven to be a very useful method of handcrafting admissible extensions for certain types of Bartnik data. We start by presenting the general strategy to obtain Bartnik mass estimates (using a Mantoulidis--Schoen type construction).

Let $(\S \cong \bS^2,g,{\mathcal{D}})$ be given 2-dimensional Bartnik data, where $g\in\M$ (some suitable set of metrics satisfying a curvature constraint) and ${\mathcal{D}}$ is a tuple formed by the remaining quantities to be prescribed at the boundary. For example, for minimal Bartnik data ${\mathcal{D}}=H\equiv 0$, while for minimal charged Bartnik data ${\mathcal{D}}=(H\equiv 0,Q)$. Let $\mathcal{A}$ denote the class of admissible extensions, which depends on the setting we are interested in, and let $\m_\textnormal{H}^{\mathcal{A}}$ and $\m_\textnormal{B}^{\mathcal{A}}$ denote the corresponding notions of Hawking and Bartnik mass, see Section~\ref{sec-intro}. 

\noindent The following procedure allows to obtain an upper bound for $\m_\textnormal{B}^{\mathcal{A}}(\S \cong \bS^2,g,{\mathcal{D}})$:
\begin{enumerate}[(I)]
\item Prove the existence of an area-preserving smooth path of metrics in $\M$ that connects $g$ to a round metric.
\item Construct a suitable collar extension of $(\S \cong \bS^2,g,{\mathcal{D}})$ in such a way that ${\mathcal{D}}$ is prescribed appropriately at the inner boundary $\S_0$; additionally arrange the inner boundary $\S_0$ to be outer-minimizing. Estimate the corresponding Hawking mass (which depends on the type of extensions we are considering), $\m_\textnormal{H}^{\mathcal{A}}$, at the outer boundary $\S_1$ of the collar extension.
\item Use a gluing tool to connect the collar extension to a model space of total mass $m > \m_\textnormal{H}^{\mathcal{A}}(\S_1)$. Since the resulting manifold is an admissible extension, we obtain
\begin{align*}
\m_\textnormal{B}^{\mathcal{A}}(\S \cong \bS^2,g,{\mathcal{D}}) \leq m.
\end{align*}
\end{enumerate}

\medskip
\section{Recent Results on Bartnik Mass Estimates}\label{sec-Bartnik-est}
In this section, we present some recent results about Bartnik mass estimates that are obtained via Mantoulidis--Schoen type constructions. For the sake of presentation, we will deviate from the historical order and first discuss the case when the Bartnik data is minimal, i.e., treat extensions of $(\S\cong \bS^2,g,H \equiv 0)$ in the asymptotically flat, asymptotically hyperbolic, and electrically charged setting. Then, we will present the results for CMC Bartnik data, i.e., for $(\S\cong \bS^2,g,H)$ with $H$ a positive constant, in the asymptotically flat and hyperbolic settings. Finally, we briefly discuss some progress regarding estimates for non-CMC Bartnik data. For a historical recollection of these results see Section~\ref{sec-MS}.

\subsection{Minimal Bartnik data}

One of the most remarkable features of the Mantoulidis--Schoen construction is that it provides a way to compute the Bartnik mass of minimal Bartnik data as was first done in~\cite{MS}. 

Motivated by the fact that horizons in a time-symmetric initial data set are stable  minimal surfaces, in~\cite{MS}, Mantoulidis and Schoen consider the set of metrics
\begin{align}\label{M+}
\M^+ &\definedas \{ g \textnormal{ is a metric on } \S \, | \, \lambda_1 ( -\Lap_{g} + K(g))> 0 \}
\end{align}
where $\lambda_1$ denotes the first eigenvalue of $ -\Lap_{g} + K(g)$, and $K(g)$ denotes the Gaussian curvature of the $2$-surface $(\S,g)$. Using the uniformization theorem, given any metric $g \in \M^+$, they obtain an area-preserving path of metrics $\{ g(t) \}_{0 \leq t \leq 1}$ in $\M^+$ connecting $g$ to $g_*$, the round metric on $\bS^2$ (see Section~\ref{sec-MS}). As discussed in Section~\ref{section-gluingmethods}, this gives rise to a family of metrics satisfying (i)-(iv) on page~\pageref{i-iv}, which we will again denote by $\{ g(t) \}_{\{0 \leq t \leq 1\}}$.

Set $v=A u(t,\cdot)$, where $A$ is a positive constant and $u>0$ is an eigenfunction corresponding to the eigenvalue $\lambda_1(t)>0$. Here, the eigenfunctions $u$ need to be chosen so that $u$ is smooth on $[0,1]\times \bS^2$ (see~\cite[Lemma A.1]{MS}). In addition, they need to be $L^2$-normalized with respect to $d\sigma_{g(t)}$. Then \eqref{eq-scal-collar} becomes
\begin{align}
\begin{split}\label{eq-scal-collar-min}
R(\g) &=\frac{2}{F(t)^2} \lambda_1(t)\\
&\qquad+ \frac{1}{A^2u^2}\(4\frac{F'(t)}{F(t)}\frac{\pr_t u}{u}  - \frac{2}{F(t)^2}\left[F'(t)^2 + 2F(t)F''(t)\right] -\a
 \).
 \end{split}
\end{align}

It is clear that, by compactness, choosing $A$ sufficiently large in terms of $u$, $F$ (and their derivatives), and $\a$, the last term of Equation~\eqref{eq-scal-collar-min} can be made arbitrarily small. Hence, we can choose $A$ so that $R(\g)>0$, since $\lambda_1(t)>0$ by assumption. However, by making $A$ larger, the length of the collar extension increases as well. In order to specify the data  near the boundary in such a way that the boundary is outer-minimizing, set $F(t)\definedas(1+\veps t^2)^{\frac{1}{2}}$, so the collar metric takes the form
\begin{align}\label{eq-collar-minimal}
\g = A_o^2u(t,\cdot)^2 dt^2+ (1+ \veps t^2)g(t),
\end{align}
where $A_o>0$ is chosen sufficiently large so that $R(\g) > 0$. Using \eqref{eq-meanc-leaf}, we know that the mean curvature at $\S_0=\{ 0 \} \times \S$ is $H(0) \equiv 0$, and positive for $\S_t$ with $t>0$. One can check directly that by taking $\veps$ small, the area $|\S_1|_{g(1)}$ can be made arbitrarily close to $|\S|_{g_o}$. 

Using property (ii) of $\{ g(t) \}_{0 \leq t \leq 1}$ from page~\pageref{i-iv}, and since $g(1)$ is a round metric, one can directly check that $u \defeq u(t,\cdot)$ is a constant for $t \in [\theta,1]$. Then one can perform the change of variables $s = A_o u t$ to write the metric $\g$ in the form required by Proposition~\ref{prop-gluing-CCMM}. Therefore, it follows that for any $m > \sqrt{\frac{|\S|_{g_o}}{16 \pi}}$, there is an admissible extension of $(\S \cong \bS^2,g_o,H_o \equiv 0)$ with ADM mass equal to $m$. Since in this case the Riemannian Penrose Inequality \eqref{eq-RPI-bound} establishes that $\sqrt{\frac{|\S|_{g_o}}{16 \pi}}$ is a lower bound for the Bartnik mass of the Bartnik data $(\S \cong \bS^2,g_o,H_o \equiv 0)$, we immediately obtain the following theorem.

\begin{theorem}[Bartnik mass of minimal Bartnik data (c.f.~{\cite[Theorem 2.1]{MS}})]
Let $(\S \cong \bS^2,g_o,H_o \equiv 0)$ be minimal Bartnik data satisfying $\lambda_1 > 0$, where $\lambda_1$ is the first eigenvalue of $-\Lap_{g_o} + K(g_o)$, with $K(g_{o})$ denoting the Gaussian curvature of $(\S,g_{o})$. Then
\begin{align*}
\m_\textnormal{B}(\S \cong \bS^2,g_o,H_o \equiv 0) = \sqrt{\frac{|\S|_{g_o}}{16 \pi}}.
\end{align*}
\end{theorem}

In the asymptotically hyperbolic setting {(see Section~\ref{subsec-intro-AH}), we are interested in constructing asymptotically hyperbolic extensions with scalar curvature bounded below by $-6$. 
In~\cite{CCM}, McCormick and the authors considered the set $\M^+$ above and the set 
\begin{align*}
\mathscr{K}^{-}&\definedas\{ g \textnormal{ is a metric on }  \S \, | \, K(g)> -3 \}.
\end{align*}
The existence of the smooth area-preserving path in $\mathscr{K}^{-}$ connecting $g$ to a round metric is given by Ricci flow in dimension two studied by Hamilton~\cite{HamiltonRF}. Again, denote by $\{ g(t) \}_{0 \leq t \leq 1}$ the modified family of metrics satisfying (i)-(iv) from page~\pageref{i-iv}.

Using the same type of collar extensions given by \eqref{eq-collar-minimal}, the condition $R(\g) \geq -6$ is equivalent via \eqref{eq-scal-collar-min} to
\begin{align*} 
A^2u^2(\lambda_1(t) + 3F(t)^2) + 2\frac{F'(t)}{F(t)}\frac{\pr_t u}{u} - \frac{1}{F(t)^2}\left[F'(t)^2 + 2F(t)F''(t)\right] -\frac{1}{2}\a \geq 0.
\end{align*}
Since $F(t) \geq 1$ for all $t \in [0,1]$, it is clear that by requiring $\{ g(t) \}_{0\leq t \leq 1}$ to lie in $\M^+$ or $\mathscr{K}^-$, it suffices to make $R(\g) > -6$ (by taking $A$ sufficiently large). Similarly as in the asymptotically flat case, these collar extensions together with Proposition~\ref{prop-gluing-CCM} produce admissible extensions, ultimately leading to the following estimate.

\begin{theorem}[Asymptotically hyperbolic Bartnik mass estimate for minimal Bartnik data{~\cite[Theorem 1.1]{CCM}}]
Let $(\S \cong \bS^2,g_o,H_o \equiv 0)$ be minimal Bartnik data satisfying either $\lambda_1 > 0$, where $\lambda_1$ is the first eigenvalue of the operator $-\Lap_{g_o}+K(g_o)>0$, or $K(g_o)>-3$, recalling that $K(g_{o})$ denotes the Gaussian curvature of $(\S,g_{o})$. Then we find the estimate
\begin{align*}
\m_\textnormal{B}^{\AH} (\S \cong \bS^2,g_o,H_o \equiv 0) \leq  \sqrt{ \frac{|\S|_{g_o}}{16\pi}}\(1 + \frac{|\S|_{g_o}}{4 \pi}  \).
\end{align*}
\end{theorem}

If, in addition, we suppose that the asymptotically hyperbolic Riemannian Penrose Inequality \eqref{eq-RPIAH-bound} holds, we would obtain the exact value
\begin{align*}
\m_\textnormal{B}^{\AH} (\S \cong \bS^2,g_o,H_o \equiv 0) =  \sqrt{ \frac{|\S|_{g_o}}{16\pi}}\(1 + \frac{|\S|_{g_o}}{4 \pi}  \).
\end{align*}
\medskip

For the Einstein--Maxwell setting, an electric charge $Q$ is part of the given Bartnik data, see Section~\ref{subsec-intro-EM}. Given charged Bartnik data $(\S \cong \bS^2,g,H \equiv 0,Q)$ where $Q$ is a real number and $g \in \M^+$ defined as in \eqref{M+}, Alaee and the authors~\cite{ACC} consider the path of metrics $\{ g(t) \}_{0 \leq t \leq 1}$ given by the uniformization theorem and satisfying (i)-(iv)  from page~\pageref{i-iv} as before. Additionally, we define the useful quantity
\begin{equation} \label{eq-kappa}
\kappa \defeq \inf_{\S \times [0,1]} \lambda_1 (t),
\end{equation}
where $\lambda_1(t) \defeq \lambda_1(-\Lap_{g(t)} + K(g(t)))$ denotes the first eigenvalue of the operator $-\Lap_{g(t)} + K(g(t))$ which smoothly depends on $t\in[0,1]$.

In this case, the main difference to the uncharged case treated in~\cite{MS} is that in addition to the collar extension \eqref{eq-collar-minimal}, it is necessary to construct a suitable divergence free electric vector field $E$ along the collar, such that $R(\g) \geq 2|E|^2_{\g}$. Using \eqref{eq-scal-collar-min} and taking $A$ sufficiently big as before, we can arrange the collar metric to satisfy  $R(\g) > 2|E|^2_{\g}$ and simultaneously make $m_{\textnormal{H}}^{\textnormal{CH}}(\S_1)$ arbitrarily close to the optimal value in  \eqref{eq-RPICH-bound} for a suitable choice of $E$, which together with Proposition~\ref{prop-gluing-ACC} gives the following result.

\begin{theorem}[Charged Bartnik mass estimate for minimal charged Bartnik data (c.f. {~\cite[Theorem~5.1]{ACC}})]
Let $(\S \cong \bS^2,g_o,H_o \equiv 0,Q_o)$ be charged Bartnik data satisfying $\lambda_1 > 0$, where $\lambda_1$ is the first eigenvalue of $-\Lap_{g_o} + K(g_o)$ with $K(g_{o})$ denoting the Gaussian curvature of $(\S,g_{o})$, $Q_o < r_o^2$, and  $\kappa > \frac{Q_o^2}{r_o^4}$, where $\kappa$ is given in \eqref{eq-kappa}. Then,
\begin{align*}
\m_\textnormal{B}^{\CH}(\S \cong \bS^2,g_o,H_o \equiv 0,Q_o) = \sqrt{\frac{|\S|_{g_o}}{16 \pi}} + \sqrt{\frac{\pi}{|\S|_{g_o}}}Q_o^2.
\end{align*}
\end{theorem}

As we have seen, for minimal Bartnik data, Mantoulidis--Schoen type extensions can be arranged to arbitrarily approach the optimal value in the corresponding Riemannian Penrose Inequality (or rather Conjecture in the asymptotically hyperbolic case). This is partly due to the fact that the size of the constant $A$ does not play a role in enforcing zero mean curvature at the inner boundary (see \eqref{eq-meanc-leaf}), which we will see below to be rather different in the case of CMC Bartnik data.

\subsection{CMC Bartnik data}\label{subsec-CMC} In this part, we restrict our attention to CMC Bartnik data $(\S \cong \bS^2,g,H)$, that is, we assume that $H\equiv H_{o}$ is a positive constant. In addition, we assume that $K(g) > 0$. We remind the reader that outer-minimizing Bartnik data $(\S\cong\bS^{2},g,H)$ satisfy the generalized Penrose inequality~\eqref{genPI}, namely
\begin{equation*}
\m_\textnormal{H}(\S \cong \bS^2,g,H) \leq \m_\textnormal{B}(\S \cong \bS^2,g,H).
\end{equation*}

Now, given CMC Bartnik data $(\S \cong \bS^2,g_o,H_o)$, the set of metrics considered by McCormick, Miao, and the authors in~\cite{CCMM} is
\begin{align*}
\mathscr{K}^{+}&\definedas\{ g \textnormal{ is a metric on }  \S \, | \, K(g)> 0 \}.
\end{align*}
As before, the uniformization theorem is used to show the existence of the desired area-preserving path, and we denote by $\{ g(t) \}_{0 \leq t \leq 1}$ the modified path satisfying (i)-(iv) from page~\pageref{i-iv}. 
For the collar metric, set $v\definedas A$, where $A$ is a positive constant, so that
\begin{align*}
\g = A^2 dt^2 + F(t)^2 g(t).
\end{align*}
Then \eqref{eq-scal-collar} takes the form
\begin{align}\label{eq-scal-collar-cmc}
R(\g) &=  \frac{2}{F(t)^2}  K(g(t)) + \frac{1}{A^2}\(  - \frac{2}{F(t)^2}\left[F'(t)^2 + 2F(t)F''(t)\right] -\a \).
\end{align}
As before, by making $A$ large in terms of $F$ (and its derivatives) and $\a$, the last term can be made arbitrarily small. In contrast to the minimal case (where the  size of $A$ played no role in prescribing the mean curvature at the boundary), the size of $A$ will now affect the  size of the mean curvature at the boundary, since by~\eqref{eq-meanc-leaf}, we have
\begin{align*}
H(0)= \frac{2}{ A}F'(0).
\end{align*}
Thus, by making $A$ large, we necessarily restrict the freedom of $H_o$. The goal is then to produce collar extensions in which $A$ is chosen in some optimal sense, and such that they appropriately propagate the Hawking mass of the given data along the collar (instead of the area which needed to be propagated appropriately in the minimal case). Using the collar extensions developed by Miao and Xie~\cite{MX}, the Hawking mass at the end of the collar can be estimated in terms of the Hawking mass of $(\S,g_o,H_o)$, $\alpha$ and $\beta$, see Remark~\ref{rem:ab}. The idea in~\cite{MX} is to use a piece of the neck of a suitable Schwarzschild manifold as the collar extension.

More precisely, let $m \in (-\infty,\frac{r_o}{2}]$ where $r_o\definedas\sqrt{ \frac{|\S|_{g_o}}{4\pi}}$ denotes the area radius of the Bartnik data. The metric of part of the neck part of the (potentially negative mass) Schwarzschild manifold corresponding to $[r_o,\infty) \times \bS^2$ can be written as before as
\begin{align*}
\g_m = ds^2 + u_{m,r_o}(s)^2 \, g_*,
\end{align*}
where $u_{m,r_o}$ is a smooth function on $[0,\infty)$ satisfying
\begin{enumerate}[(a)]
\itemsep0.25em
\item $u_{m,r_o}(0)=r_o$,
\item $u_{m,r_o}'(s)=\sqrt{1 - \frac{2m}{u_{m,r_o}(s)}}$, and
\item $u_{m,r_o}''(s)=\frac{m}{u_{m,r_o}(s)^2}$.
\end{enumerate}

Therefore, for the collar extension
\begin{align} \label{eq-collar-CMC}
\g_c &\definedas A dt^2 + \frac{u_{m,r_o}(Akt)^2}{r_o^2} g(t),
\end{align}
with $k= \frac{H_o r_o}{2}\(1-\frac{2m}{r_o}  \)^{-1/2}$, we have by \eqref{eq-scal-collar-cmc} that
\begin{align}\label{eq-scal-CMC}
R(\g_c) \geq 2 u_{m,r_o}^{-2}\( \beta - k^2 -\frac{1}{2} \a A^{-2}u_m^2(At) \).
\end{align}
Assuming the smallness condition 
\begin{align} \label{eq-smallness}
\frac{1}{4}H_o^2 r_o^2 < \frac{\beta}{1+\a}
\end{align}
and picking $m$ so that
\begin{align*}
\frac{1}{4} H_o^2 r_o^2 < \frac{\b}{1+\a} \( 1 - \frac{2m}{r_o} \),
\end{align*}
they assert that there exists $A_o$ such that
\begin{enumerate}[(a)]
\itemsep0.25em
\item $R(\g_c) > 0$,
\item $\S_0$ is isometric to $(\S,g_o)$ with mean curvature $H_o$,
\item $\S_t$ has positive mean curvature for all $t \in [0,1]$, and
\item the Hawking mass of $\S_1$ can be estimated by
\begin{align*}
\m_\textnormal{H}(\S_1) \leq \begin{cases}
\frac{1}{4} H_or_o A_o (1-k^2) + \m_\textnormal{H}(\S,g_o,H_o),\text{ if $m<0$}, \\[.25em]
\frac{1}{2} A_ok(1-k^2) + \m_\textnormal{H}(\S,g_o,H_o),\text{ if $m \geq 0$}.
\end{cases}
\end{align*}
\end{enumerate}
All the computations above are explained in~\cite{MX}.

Recalling that for $t \in [\theta,1]$ the collar extension is rotationally symmetric, one readily sees that Proposition~\ref{prop-gluing-CCMM} is applicable for any $m_* > \m_\textnormal{H}(\S_1)$. It is possible to perform an optimality analysis on the choice of $m$ and we refer the interested reader to~\cite{MX} for details. In particular, by setting $m=0$, we obtain the following estimate.

\begin{theorem}[Bartnik mass estimate for CMC Bartnik data{~\cite[Theorem~1.1]{CCMM}}]\label{thm-CCM-est}
Let $(\S\cong \bS^2, g_o, H_o)$ be CMC Bartnik data with Gaussian curvature $K(g_o)>0$. Suppose that $\frac{1}{4}H_o^2 r_o^2 < \frac{\beta}{1+\a}$, then
\begin{align*}
\m_\textnormal{B}(\S,g_o,H_o) \leq \left[ 1 + \( \frac{\a\frac{H_o^2 r_o^2}{4}}{\b - (1+\a)\frac{H_o^2 r_o^2}{4}}\)^{1/2}  \right] \m_\textnormal{H}(\S,g_o,H_o).
\end{align*}
\end{theorem}
\medskip

In~\cite{MWX}, Miao, Y.~Wang, and Xie refined the construction of collar extensions for CMC Bartnik data by proving the existence of an optimal constant $A_o$ in \eqref{eq-collar-CMC}, which makes the right hand side of \eqref{eq-scal-CMC} equal to $0$, for a suitable choice of $m$. By using a limiting argument on $m$, they were able to obtain more general estimates for $\m_\textnormal{H}(\S_1)$ than those in~\cite{MX}, which in particular hold without the smallness assumption \eqref{eq-smallness}. They obtained the following estimate.

\begin{theorem}[Bartnik mass estimate for CMC Bartnik data (c.f.{~\cite[Theorem~1.3]{MWX}})] \label{thm-MWX-est} Let $(\S\cong \bS^2, g_o, H_o)$ be CMC Bartnik data with Gaussian curvature $K(g_o)>0$ and satisfying $\frac{H_or_o}{2} \leq 1$. Then
\begin{align*}
\m_\textnormal{B}(\S\cong \bS^2, g_o, H_o) \leq C r_o (1 + \zeta_{g_o}\frac{1}{2}H_o r_o)\zeta_{g_o}\frac{1}{2}H_o r_o + \m_\textnormal{H}(\S,g_o,H_o),
\end{align*}
where  $C$ is an absolute constant and
\begin{align*}
\zeta_{g_o} &\definedas \inf \(\frac{\a}{2\b}\)^{1/2},
\end{align*}
where the infimum is taken over the set of paths $\{ g(t) \}_{0\leq t\leq 1}$ in $\mathscr{K}^+$ connecting $g_o$ to a round metric and satisfying $\tr_{g(t)}g'(t) \equiv 0$ for $t \in [0,1]$.
\end{theorem}

\begin{remark}
When the value of $H_o$ is small, the estimate in Theorem~\ref{thm-MWX-est} improves that in Theorem~\ref{thm-CCM-est}; for a detailed comparison of the estimates, see~\cite{MWX}.
\end{remark}

\begin{remark}
In~\cite{LS}, Lin and Sormani obtained interesting estimates for CMC Bartnik data $(\S\cong \bS^2,g_o,H_o)$ with normalized area $|\S_o|_{g_{o}}=4\pi$. Even though their method also relies on the construction of asymptotically flat extensions, their extensions and the character of their estimates are different in nature to the ones discussed here. For a precise statement we refer the reader to~\cite{LS}.
\end{remark}

In the asymptotically hyperbolic case, by using the collar extensions constructed by Miao and Xie~\cite{MX} and radial profile extensions that mimic a piece of an AdS-Schwarz\-schild manifold neck, the following estimate is obtained in~\cite{CCM}.

\begin{theorem}[Asymptotically hyperbolic Bartnik mass estimate for CMC Bartnik data{~\cite[Theorem~1.3]{CCM}}]  Let $(\S\cong \bS^2, g_o, H_o)$ be CMC Bartnik data with Gaussian curvature $K(g_o)>-3$. Assume further that their Hawking mass satisfies
\begin{align*}
\m_\textnormal{H}^{\AH} (\Sigma,g_o,H_o) > - \( \frac{|\S|_{g_o}}{4 \pi}  \)^{\frac{3}{2}}.
\end{align*}
Then its asymptotically hyperbolic Bartnik mass satisfies
\begin{align*}
\m_\textnormal{B}^{\AH}(\Sigma,g_o,H_o) \leq \m_\textnormal{H}^{\AH} (\Sigma,g_o,H_o) + \mathcal{E}_1(H_o,\a),
\end{align*}
where $ \mathcal{E}_{1}(H_o,\a) \to 0$ as $H_o$ or $\a$ tend to $0$, with $\a$ given by \eqref{eq-alpha}.
\end{theorem}

Using a similar approach as in the asymptotically flat case, Miao, Y.~Wang, and Xie~\cite{MWX} also refined the construction of the collar extensions in this setting and correspondingly obtained an estimate for the asymptotically hyperbolic Bartnik mass of CMC Bartnik data under the condition that the Bartnik data bound a Riemannian domain with a negative lower bound on the scalar curvature. Even though the estimate in~\cite{MWX} holds for asymptotically hyperbolic extensions with scalar curvature bounded below by $-6\kappa$, for a parameter $\kappa >0$, for the sake of exposition we state it for $\kappa=1$.

\begin{theorem}[Asymptotically hyperbolic Bartnik mass estimate for CMC Bartnik data (c.f.~{\cite[Theorem~4.1]{MWX}})] 
Let $(\S\cong \bS^2, g_o, H_o)$ be CMC Bartnik data with $K(g_o)>-3$ bounding a compact domain $\Omega$ with scalar curvature greater than or equal to $-6$. Then its asymptotically hyperbolic Bartnik mass satisfies
\begin{align*}
\m_\textnormal{B}^{\AH}(\Sigma,g_o,H_o) \leq \m_\textnormal{H}^{\AH} (\Sigma,g_o,H_o) + \mathcal{E}_2(H_o,\xi_o),
\end{align*}
where $ \mathcal{E}_{2}(H_o,\xi_o) \to 0$ as $H_o$ tends to $0$. Here, $\xi_o$ is a constant that depends only on $g_o$.
\end{theorem}

To see the exact form of $\mathcal{E}_1$ and $\mathcal{E}_2$, as well as the definition of $\xi_o$, we invite the reader to consult~\cite{CCM} and~\cite{MWX}, respectively.

\subsection{Non-CMC Bartnik data}
Suppose that $(\S \cong \bS^2,g,H)$ are given Bartnik data for which $H \geq 0$ is not necessarily constant. What can we say about estimates of their Bartnik mass? Being able to construct explicit asymptotically flat extensions with the techniques discussed here seems to be a challenging problem, and it is not clear whether it is possible at all. However, it is possible to use the Bartnik mass estimates for CMC Bartnik data described above to obtain estimates for non-CMC Bartnik data, as was pursued by McCormick~\cite{McCormick-B}.

As was mentioned in Section~\ref{sec-intro}, there are variants of the definition of Bartnik mass in the literature. These differences are mainly related to the boundary conditions of the extensions and to what is commonly called a ``non-degeneracy" condition imposed on the set of admissible extensions (see~\cite{Jauregui-B,McCormick-B} and the references cited therein). It is thus of high interest to reconcile these definitions, which is a non-trivial problem. This reconciliation has been done for at least two cases, independently by Jauregui~\cite{Jauregui-B} and McCormick~\cite{McCormick-B}. In one case, studied by McCormick~\cite{McCormick-B}, a convexity condition is imposed on Bartnik data $(\S\cong \bS^2,g,H>0)$, which is given by (c.f. \cite{McM})
\begin{equation} \label{eq-conv}
2K(g) -2H \Lap_{g}\left( \frac{1}{H} \right) - \frac{1}{2} H^2 > 0.
\end{equation}
The other case, considered by Jauregui~\cite{Jauregui-B}, is when the ``non-degeneracy" condition consists of assuming the boundary of the asymptotically flat extensions to be \emph{strictly} outer-minimizing --- note that in this survey we only assume the extensions to be (weakly) outer-minimizing. More precisely, given Bartnik data $(\S \cong \bS^2,g,H>0)$, we say that they are \emph{locally extendable} if it is possible to construct a collar extension of the Bartnik data, with non-negative scalar curvature as in Section~\ref{section-gluingmethods}.

Under the corresponding assumptions, Jauregui and McCormick show that different definitions of Bartnik mass coincide. McCormick also obtains estimates for the Bartnik mass of non-CMC Bartnik data \cite{McCormick-B}. Combining their results, we have the following estimate.
\begin{theorem}[Bartnik mass estimate for non-CMC data (cf.~{\cite{McCormick-B,Jauregui-B}})] \label{thm-nonCMC} Let $(\S \cong \bS^2,g,H)$ be Bartnik data with $H>0$. Assume either that the Bartnik data $(\S \cong \bS^2,g,H)$ are locally extendable or that the convexity condition \eqref{eq-conv} holds. Then
\begin{equation}
\m_\textnormal{B}(\S,g,H) \leq \m_\textnormal{B}(\S,g,\min_{\S} H).
\end{equation}
\end{theorem}

\begin{remark}
In Theorem \ref{thm-nonCMC}, it is implicitly understood that if the Bartnik data $(\S \cong \bS^2,g,H)$ are assumed to be locally extendable, then the extra condition that $\pr M$ be strictly outer-minimizing in the definition of admissible extensions and thus of the Bartnik mass \eqref{eq-Bartnik-mass} is required. If the convexity condition \eqref{eq-conv} is assumed, then the estimate works for the definition of Bartnik mass in \eqref{eq-Bartnik-mass} used in this survey.
\end{remark}

\begin{remark}
In~\cite{WD}, using a different methodology to the ones discussed here, Wiygul estimated the Bartnik mass for almost CMC Bartnik data $C^2$-close to the standard sphere. His techniques allowed him to study the asymptotic behavior of the Bartnik mass for small spheres. See~\cite{WD} for details.
\end{remark}

\section{Conclusions and Open Problems} \label{sec-conclusions}

The Mantoulidis--Schoen construction~\cite{MS} is a novel way of handcrafting extensions of 2-dimensional Riemannian manifolds $(\S \cong \bS^2,g)$, suited to produce admissible extensions in the context of Bartnik mass. Moreover, it is optimal in the sense that it eventually leads to a lower bound on the Bartnik mass of minimal Bartnik data $(\S\cong \bS^2,g,H\equiv 0)$, which, together with the Riemannian Penrose Inequality, establishes that 
\begin{align*}
\m_\textnormal{B}(\S,g,H\equiv 0)=\sqrt{\frac{|\S|_{g}}{16\pi}},
\end{align*}
provided that $\lambda_1(-\Lap_g + K(g))>0$. In particular, it also suggests that the Riemannian Penrose Inequality is unstable, in the sense that a manifold can almost achieve equality in \eqref{eq-RPI-bound}, while being far away from being a Schwarzschild manifold (since the minimal Bartnik data, i.e., the horizon, can be arranged to be highly non-round).

We have seen that this procedure extends to higher dimensions~\cite{CM}, and remains optimal in the context of the higher dimensional Riemannian Penrose Inequality, thus suggesting the instability of the higher dimensional Riemannian Penrose Inequality in the same sense as above. It is interesting to know that this technique led to similar instability phenomena in the the asymptotically hyperbolic \cite{G} and the electrically charged \cite{P} cases.%It will be interesting to use this fact to test whether these extensions lead to similar instability phenomena in other settings (for example the asymptotically hyperbolic and the electrically charged cases discussed above). 

Currently, the authors are using the gluing tools developed in~\cite{CCMM} to produce minimizing sequences for the Bartnik mass in the minimal case~\cite{CC}. This will shed light on subtle aspects of the stability of the Riemannian Penrose Inequality and will likely have implications on our understanding of the role of the outer-minimizing condition in the definition of Bartnik mass. Similar considerations should apply to higher dimensions. 

We have also seen that the Mantoulidis--Schoen construction can be adjusted to compute the Bartnik mass and to give Bartnik mass estimates for minimal Bartnik data in the electrically charged and asymptotically hyperbolic settings, respectively.

Modifications of the Mantoulidis--Schoen construction allow to obtain estimates for the Bartnik mass of CMC Bartnik data~\cite{CCMM}. However, the estimates obtained there are not optimal in any sense. In the light of the gluing tools developed in~\cite{CCMM}, is plausible that by picking a different type of collar extensions, the Bartnik mass estimates for CMC Bartnik data could be improved. It is evident from the work of Miao and Xie~\cite{MX} and the recent work of Miao, Y.~Wang, and Xie~\cite{MWX}, that constructing (and improving) collar extensions is a very challenging problem.

As discussed briefly in this survey, two recent works by Jauregui~\cite{Jauregui-B} and McCormick~\cite{McCormick-B} also address questions related to the outer-minimizing condition in the definition of Bartnik mass.  Remarkably, in both works, ideas \`a la Mantoulidis--Schoen are used.

We would like to point out that there have been other applications of the ideas by Mantoulidis and Schoen and their modifications different from those discussed in this survey: the works of Anderson and Jauregui~\cite{AJ}, Li and Mantoulidis \cite{LM}, Mantoulidis and Miao \cite{MM-d,MM}, Mantoulidis, Miao and Tam \cite{MMT}, and McCormick and Miao \cite{McM}. Unfortunately, discussing those in detail would go far beyond the scope of this survey.

All the considerations in this survey, and indeed all the relevant works known to the authors using, generalizing, or applying the Mantoulidis--Schoen construction, focus on time-symmetric initial data sets for (or time-slices of) the Einstein Equations --- possibly in higher dimensions, or assuming different asymptotics related to a cosmological constant, or incorporating (electric) charge. This leads to two natural questions:

\subsection*{Going beyond time-symmetry} It would be natural to study whether a similar construction of initial data sets satisfying the dominant energy condition can be performed without assuming time-symmetry. Ideally, such a construction would give us a better understanding of the (conjectured) Penrose Inequality for initial data sets beyond the Riemannian (i.e., time-symmetric) context. Furthermore, it could allow to compute admissible extensions for a non-time-symmetric version of Bartnik mass~\cite{Bartnik-mass} and thus potentially allow to compute or estimate this mass in special cases. In the $3$-dimensional, asymptotically flat, non-charged setting, this would entail constructing Riemannian manifolds $(M^{3},\g)$ carrying a symmetric $(0,2)$-tensor field $K$ to be thought of as the second fundamental form of $(M,\g)$ in the ambient spacetime. The dominant energy condition would then read
\begin{align}\label{DECK}
R(\g)-|K|_{\g}^{2}+(\tr_{\g}K)^{2}&\geq |\diver_{\g}K-d(\tr_{\g}K)|_{\g},
\end{align}
and the definition of asymptotic flatness would extend to $K$ by requesting that $K_{ij}=\mathcal{O}_{1}(r^{-p-1})$ as $r\to\infty$ while exchanging integrability of $R(\g)$ for integrability of both sides of~\eqref{DECK}. The corresponding asymptotic model space would be the so-called Kerr initial data sets, arising as special time-slices of the important (sub-extremal) Kerr spacetimes modelling rotating black holes or celestial bodies of prescribed mass $m$ and angular momentum $a$ with $m\geq |a|$, in vacuum, see e.g. Wald~\cite{Wald}. These are not spherically, but still axially symmetric. Minimal Bartnik data would need to be replaced by ``MOTS'' Bartnik data, for which in particular $H=\tr_{\S}K$ would be requested to hold on $\S\cong\bS^{2}$.
Constructing admissible extensions for MOTS Bartnik data would thus require
\begin{enumerate}[(a)]
\itemsep0.25em
\item to come up with a new idea of how to construct suitable axisymmetric ``collars'', including a construction for $K$, near MOTS, replacing the condition $R(\g)>0$ by the (strict) version of~\eqref{DECK}, and
\item to be able to smoothly glue such a collar to a suitable time-slice of a sub-extremal Kerr spacetime with suitably chosen mass $m$, while ensuring~\eqref{DECK} along the gluing bridge and that the inner boundary of the constructed initial data set is outer-minimizing.
\end{enumerate}

\subsection*{Studying the time-evolution of the constructed initial data} From a physics perspective, it is natural to ask what spacetimes arising by time-evolution under the Einstein Equations from initial data sets such as those described in this survey will look like. This question silently assumes that one chooses what is a called a \emph{matter model} describing the matter present in the gravitating system one wishes to model mathematically. Mantoulidis and Schoen~\cite{MS} as well as Alaee and the authors~\cite{ACC} give examples of matter models which are compatible with time-evolution of the respective constructed time-symmetric initial data sets. It would be desirable to know more generally which matter models are compatible with time-evolution under the Einstein Equations for the various generalizations of the Mantoulidis--Schoen construction described in this survey.

\bibliographystyle{amsplain}
\bibliography{bib-prelim}

\end{document}